%

%
%
%
%
\documentstyle{amsppt}
\loadbold
\def\cstar{$C^*$-algebra}

\def\<{\left<}										
\def\>{\right>}

\def\tr{\text{trace}\,}
\def\cp{\text{cp}}
\def\seq{\text{seq}}


\magnification=\magstephalf

\topmatter
\title
Subalgebras of $C^*$-algebras III:\\
Multivariable operator theory
\endtitle

\author William Arveson
\endauthor

\affil Department of Mathematics\\
University of California\\Berkeley CA 94720, USA
\endaffil

\date 18 May 1997
\enddate
\thanks This research was supported by
NSF grant DMS-9500291
\endthanks
%
%
%
\abstract 
A $d$-contraction is a $d$-tuple $(T_1,\dots,T_d)$ 
of mutually commuting operators acting on a common 
Hilbert space $H$ such that 
$$
\|T_1\xi_1+T_2\xi_2+\dots +T_d\xi_d\|^2\leq 
\|\xi_1\|^2+\|\xi_2\|^2+\dots+\|\xi_d\|^2
$$
for all $\xi_1,\xi_2,\dots,\xi_d\in H$.  
These are the higher dimensional counterparts of
contractions.  We show that 
many of the operator-theoretic aspects of function 
theory in the unit disk generalize to 
the unit ball $B_d$ in complex $d$-space, 
including von Neumann's inequality and the
model theory of contractions.  
These results depend on properties 
of the $d$-shift, a distinguished $d$-contraction which 
acts on a new $H^2$ space associated with $B_d$, and 
which is the higher dimensional counterpart 
of the unilateral shift.   $H^2$ and the $d$-shift are 
highly unique.  Indeed, by exploiting 
the noncommutative Choquet boundary 
of the $d$-shift relative to its generated
$C^*$-algebra we find that
there is more uniqueness in dimension $d\geq 2$ 
than there is in dimension one.  
\endabstract


\toc
\specialhead{Introduction} 
\endspecialhead

\specialhead{Part I: Function theory}
\endspecialhead

\subhead 1.  Basic properties of $H^2$
\endsubhead
\subhead 2.  Multipliers and the $d$-dimensional shift
\endsubhead
\subhead 3.  von Neumann's inequality and the sup norm
\endsubhead
\subhead 4.  Maximality of the $H^2$ norm
\endsubhead

\specialhead{Part II: Operator theory}
\endspecialhead

\subhead 5.  The Toeplitz $C^*$-algebra
\endsubhead
\subhead 6.  $d$-contractions and $A$-morphisms
\endsubhead
\subhead 7.  The $d$-shift as an operator space
\endsubhead
\subhead 8.  Various applications
\endsubhead
\specialhead{Appendix: Trace estimates}
\endspecialhead
\specialhead{References} 
\endspecialhead
\endtoc

\endtopmatter
%

\document

\subheading{Introduction}

This paper concerns function theory and operator theory 
relative to the unit ball in complex $d$-space 
$\Bbb C^d$, $d=1,2,\dots$.  A $d$-contraction is a 
$d$-tuple $(T_1,\dots,T_d)$ of mutually commuting operators 
acting on a common Hilbert space $H$ satisfying 
$$
\|T_1\xi_1+\dots+T_d\xi_d\|^2\leq \|\xi_1\|^2+\dots+\|\xi_d\|^2, 
$$
for every $\xi_1,\dots,\xi_d\in H$.  This inequality 
simply means that
the ``row operator" defined by the $d$-tuple, viewed as 
an operator from the direct sum of $d$ copies of $H$ to $H$, 
is a contraction.  It is essential that the component 
operators commute with one another.  

We show that 
there exist $d$-contractions which are not polynomially
bounded in the sense that there is no constant $K$ satisfying 
$$
\|f(T_1,\dots,T_d)\|\leq K
\sup\{|f(z_1,\dots,z_d)| : |z_1|^2+\dots+|z_d|^2\leq 1\}
$$
for every polynomial $f$.  In fact, we single 
out a particular $d$-contraction $(S_1,\dots,S_d)$ (called 
the $d$-shift) which is not polynomially 
bounded but which gives rise to the appropriate 
version of von Neumann's inequality 
with constant $1$: for every $d$-contraction $(T_1,\dots,T_d)$
one has  
$$
\|f(T_1,\dots,T_d)\|\leq\|f(S_1,\dots,S_d)\|
$$
for every polynomial $f$.  Indeed the indicated homomorphism 
of commutative operator algebras is completely contractive.  

The $d$-shift acts naturally on a space of holomorphic 
functions defined on the 
open unit ball $B_d\subseteq \Bbb C^d$, which we call $H^2$.  
This space is a natural generalization of the familiar Hardy
space of the unit disk, but it differs from other 
``$H^2$" spaces in several  
ways.  For example, unlike the space $H^2(\partial B_d)$ 
associated with normalized surface area on the sphere or the
space $H^2(B_d)$ associated with volume measure over the 
interior, $H^2$ is not associated with any measure on 
$\Bbb C^d$.  Consequently, the associated multiplication operators 
(the component operators of the $d$-shift) do not form 
a subnormal $d$-tuple.  Indeed, since the naive 
form of von Neumann's inequality described 
above fails, no effective model 
theory in dimension $d\geq 2$ could be based 
on subnormal operators.  Thus by giving up 
the requirement of subnormality for models, one gains a 
theory in which models not only exist in all 
dimensions but are unique as well.  

In the first part of this paper we work out the basic theory
of $H^2$ and its associated multiplier algebra, and we show
that the $H^2$ norm is the largest Hilbert 
norm on the space of polynomials which is appropriate for 
the operator theory of $d$-contractions.  

In part II we 
emphasize the role of ``$\Cal A$-morphisms".  These are 
completely positive linear maps of the $d$-dimensional 
counterpart of the Toeplitz \cstar\ which bear a particular
relation to the $d$-shift.  Every 
$d$-contraction corresponds to a {\it unique} $\Cal A$-morphism, 
and on that observation we base a model theory for 
$d$-contractions which provides an appropriate generalization
of the Sz.-Nagy Foias theory of contractions  \cite{43}
to arbitrary dimension $d\geq 1$ (see section 8).  In section 7 we 
introduce a sequence of numerical invariants $E_n(\Cal S)$,
$n=1,2,\dots$ for arbitrary operator spaces $\Cal S$.  
We show that the $d$-dimensional operator space $\Cal S_d$ 
generated by the $d$-shift is maximal in the sense that 
$E_n(\Cal S_d)\geq E_n(\Cal S)$ for every $n\geq 1$ and
for every $d$-dimensional
operator space $\Cal S$ consisting of mutually commuting 
operators.  More significantly, we show that when
$d\geq 2$, $\Cal S_d$ is {\it characterized} by
this maximality property.  That characterization
fails for single operators 
(i.e. one-dimensional operator spaces).  
We may conclude that, perhaps contrary to one's 
function-theoretic intuition, there is more uniqueness 
in dimension $d\geq 2$ than there is in dimension one.

Since this paper is a logical sequel to \cite{3}, \cite{4} and 
so many years have passed since the publication of 
its two predecessors, it seems appropriate to comment on 
its relationship to them.  On the one hand, we have 
come to the opinion that the program proposed in \cite{4, chapter 1} 
for carrying out dilation theory in higher dimensions 
must be modified.  That program 
gives necessary and sufficient conditions for finding 
normal dilations in multivariable operator theory.  
However, the results below provide
two reasons why normal dilations are inappropriate for 
commutative sets of 
operators associated with the unit ball $B_d$.  
First, they may not exist (a $d$-contraction 
need not have a normal dilation
with spectrum in $\partial B_d$, c.f. Remark 3.13)
and second, when they do exist they 
are not unique (there can be
many normal dilations of a given $d$-contraction 
which have the stated properties but 
which are not unitarily equivalent to each other).  

On the other hand, the results of this paper also 
demonstrate that 
other aspects of the program of \cite{3,4} are  
well-suited for multivariable 
operator theory.  For example, we will see that 
boundary representations...the noncommutative 
counterparts of Choquet boundary points in the commutative 
theory of function spaces...play an important role in 
the operator theory of $B_d$.  Boundary representations
serve to explain 
the notable fact that in higher dimensions there is 
more uniqueness than there is 
in dimension one (cf. Theorem 7.7 and its 
Corollary), and they 
provide concrete information about the 
absence of inner functions for the 
$d$-shift (cf. Proposition 8.13).  

We were encouraged to return to these problems
by recent results in the theory of $E_0$-semigroups.  
There is a dilation theory by which, starting with
a semigroup of completely 
positive maps of $\Cal B(H)$, one obtains an $E_0$-semigroup
as its ``minimal dilation" \cite{14,6,7,8,9,10}.  
In its simplest form, this dilation theory starts with
a normal completely 
positive map $P: \Cal B(H)\to\Cal B(H)$ satisfying 
$P(\bold 1)=\bold 1$ and constructs from it a unique endomorphism
of $\Cal B(K)$ where $K$ is a Hilbert space containing 
$H$.  When one looks closely at this procedure 
one sees that there should be a 
corresponding dilation theory for sets of 
operators such as $d$-contractions.  

We have reported on some of these results in a conference 
at the Fields institute in Waterloo in early 1995.  That lecture
concerned the dilation theory of semigroups of completely positive 
maps, $\Cal A$-morphisms and the issue of uniqueness.  
However, at that time we had not yet reached a 
definitive formulation of the application to operator theory.  

There is a large literature relating to von Neumann's 
inequality and dilation theory for sets of operators, and no attempt
has been made to compile a comprehensive list of references here.  
More references can be found in \cite{26},\cite{27}.  
Finally, I want to thank Ra\'ul Curto for bringing 
me back up to date on the literature of multivariable 
operator theory.


\subheading{1. Basic properties of $H^2$}
Throughout this paper we will be concerned with function theory and 
operator theory as it relates to the unit ball $B_d$ in complex
$d$-dimensional space $\Bbb C^d$, $d=1,2,\dots$,
$$
B_d = \{z=(z_1,z_2,\dots,z_d)\in \Bbb C^d: \|z\|<1 \},
$$
where $\|z\|$ denotes the norm associated with the usual inner 
product in $\Bbb C^d$
$$
\|z\|^2 = |z_1|^2+|z_2|^2+\dots+|z_d|^2.  
$$

In dimension $d=1$ there is a familiar Hardy space which can be 
defined in several ways.  We begin by reiterating 
one of the definitions 
of $H^2$ in a form that we will generalize verbatim to higher 
dimensions.  Let $\Cal P$ be the algebra of all holomorphic 
polynomials $f$ in a single complex variable $z$.  Every 
$f\in\Cal P$ has a finite Taylor series expansion 
$$
f(z) = a_0 + a_1z + \dots +a_nz^n
$$
and we may define the norm $\|f\|$ of such a polynomial as the 
$\ell^2$ norm of its sequence of Taylor coefficients
$$
\|f\|^2 = |a_0|^2 + |a_1|^2 +\dots +|a_n|^2.  \tag{1.1}
$$
The norm $\|f\|$ is of course associated with an inner product
on $\Cal P$, and the completion of $\Cal P$ in this norm 
is the Hardy space $H^2$.  It is well known that
the elements of $H^2$ can 
be realized concretely as analytic functions
$$
f: \{|z|<1\}\to \Bbb C
$$
which obey certain growth conditions near the boundary
of the unit disk.  

Now consider the case of dimension $d>1$.  $\Cal P$ will denote
the algebra of all complex holomorphic polynomials $f$ in the 
variable $z=(z_1,z_2,\dots,z_d)$.  Every such polynomial 
$f$ has a unique expansion into a finite series 
$$
f(z) = f_0(z) + f_1(z) + \dots +f_n(z) \tag{1.2}
$$
where $f_k$ is a homogeneous polynomial of degree $k$.  
We refer to (1.2) as the Taylor series of $f$.  

\proclaim{Definition 1.3}
Let $V$ be a complex vector space.  By a Hilbert seminorm on 
$V$ we mean a seminorm which derives from a positive semidefinite
inner product $\<\cdot,\cdot\>$ on $V$ by way of 
$$
\|x\| = \<x,x\>^{1/2}, \qquad x\in V.  
$$
\endproclaim

We will define a Hilbert seminorm on $\Cal P$ by imitating 
formula (1.1), where $c_k$ is replaced with $f_k$.  To 
make that precise we must view the expansion (1.2) in a 
somewhat more formal way.  The space $E = \Bbb C^d$ is a 
$d$-dimensional vector space having a 
distinguished inner product 
$$
\<z,w\> = z_1\bar w_1+z_2\bar w_2+\dots+ z_d \bar w_d.   
$$
For each
$n=1,2,\dots$ we write $E^n$ for the symmetric tensor product of 
$n$ copies of $E$.  $E^0$ is defined as the one dimensional
vector space $\Bbb C$ with its usual inner product.  
For $n\geq 2$, $E^n$ is the subspace of the full tensor 
product $E^{\otimes n}$ consisting of all vectors fixed under the 
natural representation of the permutation group $S_n$,
$$
E^n = \{\xi\in E^{\otimes n}: U_{\pi}\xi = \xi, \quad\pi\in S_n\},
$$
$U_\pi$ denoting the isomorphism of $E^{\otimes n}$ defined on 
elementary tensors by 
$$
U_\pi(z_1\otimes z_2\otimes \dots\otimes z_n) = 
z_{\pi^{-1}(1)}\otimes z_{\pi^{-1}(2)}\otimes\dots\otimes z_{\pi^{-1}(n)},
\qquad z_1\in E.
$$
For a fixed vector $z\in E$ we will use the notation 
$$
z^n = z^{\otimes n}\in E^n
$$
for the $n$-fold tensor product of copies of $z$ ($z^0\in E^0$
is defined as the complex number 
$1$).  $E^n$ is linearly
spanned by the set $\{z^n: z\in E\}$, $n=0,1,2,\dots$.    

Now every homogeneous polynomial $g: E\to \Bbb C$ of degree $k$ 
determines a unique linear functional $\tilde g$ on $E^k$ by
$$
g(z) = \tilde g(z^k), \qquad z\in E
$$
(the uniqueness of $\tilde g_k$ follows
from the fact that 
$E^k$ is spanned by $\{z^k: z\in E\}$), and 
thus the Taylor series (1.2) can be written in the form 
$$
f(z) = \sum_{k=0}^n \tilde f_k(z^k), \qquad z\in E,
$$
where $\tilde f_k$ is a uniquely determined
linear functional on $E^k$ 
for each $k=0,1,\dots,n$.  
Finally, if we bring in the inner product on $E$ then $E^{\otimes k}$ 
becomes a $d^k$-dimensional complex Hilbert space.   Thus the 
subspace $E^k$ is also a finite dimensional Hilbert space in 
a natural way.  Making use of the Riesz lemma, 
we find that there is a unique vector $\xi_k\in E^k$ such
that 
$$
\tilde f_k(z^k) = \<z^k,\xi_k\>,\qquad z\in E,
$$
and finally the Taylor 
series for $f$ takes the form 
$$
f(z) = \sum_{k=0}^n \<z^k,\xi_k\>,\qquad z\in E.  \tag{1.4}
$$
We  define a Hilbert seminorm on $\Cal P$ as follows 
$$
\|f\|^2 = \|\xi_0\|^2+\|\xi_1\|^2+\dots+\|\xi_n\|^2. \tag{1.5}
$$
The seminorm $\|\cdot\|$ is obviously a norm on $\Cal P$ in 
that $\|f\|=0 \implies f=0$.  

\proclaim{Definition 1.6}
$H^2_d$ is defined as the Hilbert space obtained by completing
$\Cal P$ in the norm (1.5).  
\endproclaim

When there is no possibility of confusion concerning the
dimension we will abbreviate $H^2_d$ with the simpler $H^2$.  
We first point out that the elements of $H^2$ can be identified 
with the elements of the symmetric Fock space over $E$,
$$
\Cal F_+(E) = E^0\oplus E^1\oplus E^2\oplus \dots
$$
the sum on the right denoting the infinite direct sum of Hilbert
spaces.  

\proclaim{Proposition 1.7}
For every $f\in \Cal P$ let $Jf$ be the element of $\Cal F_+(E)$
defined by 
$$
Jf = (\xi_0,\xi_1,\dots) 
$$
where $\xi_0, \xi_1,\dots$ is the sequence of Taylor coefficients
defined in (1.4), continued so that $\xi_k=0$ for $k>n$.  
Then $J$ extends uniquely to an anti-unitary operator mapping 
$H^2$ onto $\Cal F_+(E)$.  
\endproclaim
\demo{proof}
The argument is perfectly straightforward, once one realizes 
that $J$ is not linear but anti-linear \qed
\enddemo

We can also identify the elements of $H^2$ in more concrete 
terms as analytic functions defined on the ball $B_d$: 

\proclaim{Proposition 1.8}
Every element of $H^2$ can be realized as an analytic
function in $B_d$ having a power series expansion of the form 
$$
f(z) = \sum_{k=0}^\infty \<z^k,\xi_k\> 
\qquad z=(z_1,\dots,z_d)\in B_d
$$
where the $H^2$ norm of $f$ is given by
$\|f\|^2 =\sum_k \|\xi_k\|^2 < \infty$.  Such functions 
$f$ satisfy a growth condition of the form 
$$
|f(z)| \leq \frac{\|f\|}{\sqrt{1-\|z\|^2}}, 
\qquad z\in B_d.  
$$
\endproclaim
\demo{proof}
Because of Proposition (1.7) the elements of $H^2$ can be identified
with the formal power series having the form 
$$
f(z) = \sum_{k=0}^\infty \<z^k,\xi_k\>,  \tag{1.9}
$$
where the sequence $\xi_k\in E^k$ satisfies 
$$
\sum_{k=0}^\infty \|\xi_k\|^2 = \|f\|^2 <\infty.  \tag{1.10}
$$
Because of (1.10) the series in (1.9) is easily seen to converge
in $B_d$ and satisfies the stated growth condition.  

In more detail, since the norm of a 
vector in $E^k$ of the form $z^k$, 
$z\in E$, satisfies 
$$
\|z^k\|^2 = \<z^k,z^k\>= \<z,z\>^k = \|z\|^{2k}, 
$$
we find that
$$
|\<z^k,\xi_k\>|\leq \|z^k\| \cdot\|\xi_k\| \leq \|z\|^k \|\xi_k\|,
$$
and hence for all $z\in E$ satisfying $\|z\|<1$ we have 
$$
\sum_{k=0}^\infty|\<z^k,\xi_k\>|\leq (\sum_{k=0}^\infty\|z\|^{2k})^{1/2} 
(\sum_{k=0}^\infty\|\xi_k\|^2)^{1/2} =(1-\|z\|^2)^{-1/2}\|f\|,  
$$
as asserted \qed
\enddemo

We will make frequent use of the following  
family of functions in $H^2$.  For every $x\in B_d$ 
define $u_x:B_d\to \Bbb C$ by
$$
u_x(z) = (1-\<z,x\>)^{-1}, \qquad \|z\|<1.  \tag{1.11}
$$
$u_x(z)$ is clearly analytic in $z$ and co-analytic in 
$x$.  The useful properties of the set of functions
$\{u_x: x\in B_d\}$ are 
summarized in the following proposition, which gives 
the precise sense in which $H^2$ is characterized in 
abstract terms by the positive definite reproducing 
kernel $k:B_d\times B_d\to\Bbb C$,
$$
k(x,y) = (1-\<x,y\>)^{-1}.  
$$

\proclaim{Proposition 1.12}
$u_x$ belongs to $H^2$ for every $x\in B_d$ 
and these functions satisfy
$$
\<u_x,u_y\> = (1-\<y,x\>)^{-1}.  \tag{1.13}
$$
$H^2$ is spanned by $\{u_x: x\in B_d\}$, and for every 
$f\in H^2$ we have 
$$
f(z) = \<f,u_{z}\>, \qquad z\in B_d. \tag{1.14}
$$

Moreover, if $K$ is any Hilbert space spanned by a 
subset of its elements $\{v_x: x\in B_d\}$ which satisfy
$$
\<v_x,v_y\> = (1-\<y,x\>)^{-1}, \qquad x,y\in B_d,
$$
then there is a unique unitary operator $W: H^2\to K$ 
such that $Wu_x = v_x$, $x\in B_d$.  
\endproclaim
\demo{proof}
The proof is straightforward.  For example, to 
see that $u_x$ belongs to $H^2$ we simply examine its
Taylor series
$$
u_x(z) = (1-\<z,x\>)^{-1} = \sum_{k=0}^\infty \<z,x\>^k.  
$$
Noting that $\<z,x\>^k = \<z^k,x^k\>_{E^k}$ we can 
write
$$
u_x(z) = \sum_{k=0}^\infty \<z^k,x^k\>_{E^k}.  
$$
This shows that the sequence of Taylor coefficients 
of $u_x$ is 
$$
Ju_x = (1,x,x^2,\dots)\in \Cal F_+(E).  
$$
Hence $u_x$ belongs to $H^2$ and we have 
$$
\<u_x,u_y\> = \<Ju_y,Ju_x\>_{\Cal F_+(E)} = 
\sum_{k=0}^\infty \<y,x\>^k = (1-\<y,x\>)^{-1}.  
$$
Formula (1.13) follows.  

Similarly, a direct application of Proposition
1.7 establishes (1.14).  From 
the latter it follows that $\{u_x: x\in B_d\}$ spans $H^2$.  
Indeed, if $f$ is any function in $H^2$ which is orthogonal 
to every $u_x$ then 
$$
f(z) = \<f,u_{z}\>=0 \qquad \text{for every }z\in B_d,
$$
hence $f=0$.  

Finally, the second paragraph is obvious from the fact 
that for every finite subset $x_1,\dots,x_n\in B_d$ 
and $c_1,\dots,c_n\in \Bbb C$ we have 
$$
\|c_1u_{x_1} + \dots + c_nu_{x_n}\|^2 = 
\|c_1v_{x_1}+\dots + c_nv_{x_n}\|^2, 
$$
which is apparent after expanding both sides and comparing
inner products\qed 
\enddemo

The $H^2$ norm is invariant under the natural action of 
the unitary group of $\Bbb C^d$, as summarized by

\proclaim{Corollary}
Let $V$ be a unitary operator on the Hilbert space 
$E=\Bbb C^d$.  Then there is a unique unitary operator 
$\Gamma(V)\in \Cal B(H^2)$ satisfying 
$$
\Gamma(V) u_x = u_{Vx}, \qquad x\in B_d.  \tag{1.15}
$$
$\Gamma$ is a strongly continuous unitary representation 
of $\Cal U(\Bbb C^d)$ on $H^2$ whose action on functions 
is given by 
$$
\Gamma(V)f(z) = f(V^{-1}z), \qquad z\in B_d, f\in H^2.  \tag{1.16}
$$
\endproclaim
\demo{proof}
Fix $V\in \Cal U(\Bbb C^d)$.  For any $x,y\in B_d$ we have 
$$
\<u_{Vx},u_{Vy}\> = (1-\<Vy,Vx\>)^{-1} = (1-\<y,x\>)^{-1} =
\<u_x,u_y\>.  
$$
\enddemo
It follows from Proposition 1.12 that there is a unique 
unitary operator $\Gamma(V)\in \Cal B(H^2)$ satisfying 
(1.15).  It is clear from (1.15) that 
$\Gamma(V_1V_2) = \Gamma(V_1)\Gamma(V_2)$, and strong
continuity follows from the fact that 
$$
\<\Gamma(V)u_x,u_y\> = \<u_{Vx},u_y\>=(1-\<y,Vx\>)^{-1}
$$
is continuous in $V$ for fixed $x,y\in B_d$, together with 
the fact that $H^2$ is spanned by $\{u_z: z\in B_d\}$.  

Finally, from (1.14) we see that for every $f\in H^2$
and every $z\in B_d$, 
$$
\align
f(V^{-1}z) &= \<f,u_{V^{-1}z}\> = \<f,\Gamma(V^{-1})u_z\> =
\<f,\Gamma(V)^*u_z\> \\
&= \<\Gamma(V)f,u_z\> = (\Gamma(V)f)(z),
\endalign
$$
proving (1.16)\qed


\subheading{2. Multipliers and the $d$-dimensional shift}

By a {\it multiplier} of $H^2$ we mean a complex-valued 
function $f: B_d\to \Bbb C$ with the property 
$$
f\cdot H^2\subseteq H^2.  
$$
The set of multipliers is a complex algebra of 
functions defined on the ball $B_d$ which contains 
the constant functions, and 
since $H^2$ itself contains the constant function $1$ it follows
every multiplier must belong to $H^2$.  In particular, 
multipliers are analytic functions on $B_d$.  

\proclaim{Definition 2.1}
The algebra of all multipliers is denoted $\Cal M$.  
$H^\infty$ will denote the Banach algebra of all bounded analytic
functions $f: B_d\to \Bbb C$ with norm 
$$
\|f\|_\infty = \sup_{\|z\|<1}|f(z)|.  
$$  
\endproclaim

The following result implies that $\Cal M\subseteq H^\infty$ 
and the inclusion map of $\Cal M$ in $H^\infty$ becomes a contraction
after one endows $\Cal M$ with its natural norm.  

\proclaim{Proposition 2.2}
Every $f\in \Cal M$ defines a unique bounded operator 
$M_f$ on $H^2$ by way of 
$$
M_f: g\in H^2\to f\cdot g\in H^2.  
$$
The natural norm in $\Cal M$
$$
\|f\|_\Cal M = \sup \{\|f\cdot g\|: g\in H^2, \|g\|\leq 1\}
$$
satisfies 
$$
\|f\|_\Cal M = \|M_f\|, 
$$
the right side denoting the operator norm in $\Cal B(H^2)$, 
and we have 
$$
\|f\|_\infty \leq \|f\|_\Cal M, \qquad f\in \Cal M.  
$$
\endproclaim
\demo{proof}
Fix $f\in \Cal M$.  Notice first that if $g$ is an arbitrary 
function in $H^2$ then by (1.16) we have 
$$
\<M_fg,u_z\> = \<f\cdot g,u_{z}\> = f(z)g(z).  \tag{2.3}
$$
A straightforward application of the closed graph theorem 
(which we omit) now shows that the operator 
$M_f$ is bounded.  

It is clear that $\|f\|_\Cal M = \|M_f\|$.  We claim now
that for each $x\in B_d$ one has 
$$
M_f^*u_x = \bar f(x)u_x.  \tag{2.4}
$$
Indeed, since $H^2$ is spanned by $\{u_{y}: y\in B_d\}$ it is 
enough to show that 
$$
\<M_f^*u_x, u_{y}\> = \bar f( x)\<u_x,u_{y}\>, \qquad y\in B_d.  
$$
For fixed $y$ the left side is 
$$
\<u_x,f\cdot u_y\>=\overline{\<f\cdot u_y,u_x\>}
$$ 
By (1.16) the latter is 
$$
\overline{f(x)u_y(x)}=\bar f(x)\overline{(1-\<x,y\>)^{-1}}
= \bar f(x)(1-\<y,x\>)^{-1} = \bar f(x)\<u_x,u_y\>,
$$
and (2.4) follows.  

Finally, (2.4) implies that for every $x\in B_d$ we have 
$$
|f(x)| = \|M_f^*u_x\|/\|u_x\| \leq \|M_f^*\| = \|M_f\| = \|f\|_\Cal M
$$
as required\qed
\enddemo

We turn now to the definition of the $d$-dimensional analogue
of the unilateral shift.  Let $e_1,e_2,\dots,e_d$ be an orthonormal
basis for $E=\Bbb C^d$, and define $z_1,z_2,\dots,z_d\in \Cal P$ 
by 
$$
z_k(z) = \<z,e_k\>, \qquad x\in \Bbb C^d.  
$$
Such a $d$-tuple of linear functionals will be called a 
{\it system of coordinate functions}.  If $z_1^\prime,
z_2^\prime,\dots,z_d^\prime$ is another system of coordinate 
functions then there is a unique unitary operator 
$V\in\Cal B(\Bbb C^d)$ satisfying 
$$
z_k^\prime(x) = z_k(V^{-1}x), \qquad 1\leq k\leq d, 
\quad x\in \Bbb C^d. 
\tag{2.5} 
$$

\proclaim{Proposition 2.6}
Let $z_1,z_2,\dots,z_d$ be a system of coordinate functions for
$\Bbb C^d$.  Then for every complex number $a$ and polynomials 
$f_1,f_2,\dots,f_d\in \Cal P$ we have
$$
\|a\cdot 1 + z_1f_1+\dots +z_df_d\|^2 \leq 
|a|^2 + \|f_1\|^2 +\dots +\|f_d\|^2,
$$
$\|\cdot\|$ denoting the norm in $H^2$.  
\endproclaim
\demo{proof}
We claim first that each $z_k$ is a multiplier.  Indeed, if 
$f\in H^2$ has Taylor series 
$$
f(x) = \sum_{n=0}^\infty \<z^n,\xi_n\>
$$
with $\sum_n\|\xi_n\|^2 = \|f\|^2<\infty$ then we have
$$
z_k(x)f(x) = \sum_{n=0}^\infty \<x,e_k\>\<x^n,\xi_n\>.  \tag{2.7}
$$
Now 
$$
\<x,e_k\>\<x^n,\xi_n\> = \<x^{n+1},e_k\otimes \xi_n\>.  
$$
So if $e_k\cdot \xi_n$ denotes the projection of the vector
$e_k\otimes \xi_n\in E\otimes E^n$ to the subspace $E^{n+1}$
then (2.7) becomes 
$$
z_k(x)f(x) = \sum_{n=0}^\infty \<x^{n+1},e_k\cdot \xi_n\>.  
$$
Since 
$$
\sum_{n=0}^\infty\|e_k\cdot \xi_n\|^2 \leq 
\sum_{n=0}^\infty \|\xi_n\|^2 = \|f\|^2
$$
it follows that $z_kf\in H^2$ and in fact 
$$
\|z_k f\|\leq \|f\|, \qquad f\in H^2.  
$$

Thus each multiplication operator $M_{z_k}$ is a contraction
in $\Cal B(H^2)$. Consider the Hilbert space 
$$
K = \Bbb C \oplus
\underbrace{H^2\oplus \dots\oplus H^2}_{d \text{ times}}
$$
and the operator $T: K\to H^2$ defined by 
$$
T(a,f_1,\dots,f_d) = a\cdot 1+ z_1f_1+\dots+z_df_d .  
$$
The assertion of Proposition 2.6 is that $\|T\|\leq 1$.  In 
fact, we show that the adjoint of $T$, $T^*: H^2\to K$, 
is an isometry.  A routine computation implies that for 
all $f\in H^2$ we have 
$$
T^*f = (\<f,1\>,A_1^*f,\dots,A_d^*f) \in K,
$$
where we have written $A_k$ for the multiplication operator 
$M_{z_k}$, $k=1,\dots,d$.  Hence $TT^*\in \Cal B(H^2)$ is given
by
$$
TT^* = E_0 + A_1A_1^*+\dots+ A_dA_d^*
$$
where $E_0$ is the projection on the one-dimensional space 
of all constant functions in $H^2$.  We establish the key 
assertion as a lemma for future reference.  

\proclaim{Lemma 2.8}
Let $z_1,\dots,z_d$ be a system of coordinate functions for 
$\Bbb C^d$ and let $A_k = M_{z_k}$, $k=1,2,\dots,d$.  Let 
$E_0$ be the projection onto the one-dimensional space 
of constant functions in $H^2$.  Then 
$$
E_0 + A_1 A_1^* + \dots A_d A_d^* = \bold 1.  
$$
\endproclaim
\demo{proof}
Since $H^2$ is spanned by $\{u_z: z\in B_d\}$ it is enough to 
show that for all $x,y\in B_d$ we have 
$$
\<E_0u_x,u_y\> + \sum_{k=1}^d \<A_kA_k^*u_x,u_y\> = \<u_x,u_y\>.
\tag{2.9}
$$
Since each $A_k$ is a multiplication operator, formula (2.4) 
implies that 
$$
A_k^* u_x = \bar z_k(x)u_x = \<e_k,x\>u_x,
$$
for $x\in B_d$.  Thus we can write 
$$
\align
\sum_{k=1}^d \<A_kA_k^*u_x,u_y\> &= \sum_{k=1}^d\<A_k^*u_x,A_k^*u_y\>
=\sum_{k=1}^d\<e_k,x\>\<y,e_k\>\<u_x,u_y\> \\
&=
\<y,x\>\<u_x,u_y\> = \<y,x\>(1-\<y,x\>)^{-1}.  
\endalign
$$
On the other hand, noting that $u_0 = 1$ and $\|u_0\|=1$, 
the projection 
$E_0$ is given by $E_0(f) = \<f,u_0\>u_0$, $f\in H^2$.  
Hence
$$
\<E_0u_x,u_y\> = \<u_x,u_0\>\<u_0,u_y\> = 1
$$
because $\<u_x,u_0\> = (1-\<x,0\>)^{-1} = 1$ for every 
$x\in B_d$.  It follows that 
$$
\align
\<E_0u_x,u_y\> + \sum_{k=1}^d\<A_kA_k^*u_x,u_y\> &= 
1 + \<y,x\>(1-\<y,x\>)^{-1} = (1-\<y,x\>)^{-1} \\
&= \<u_x,u_y\>, 
\endalign
$$
as asserted \qed

That completes the proof of Proposition 2.6\qed
\enddemo

\enddemo

\proclaim{Definition 2.10}
Let $z_1,\dots,z_d$ be a system of coordinate functions for 
$\Bbb C^d$ and let $S_k = M_{z_k}$, $k=1,2,\dots,d$.  
The $d$-tuple of operators 
$$
\bar S = (S_1,S_2,\dots,S_d)
$$
is called the $d$-dimensional shift or, briefly, the 
$d$-shift.  
\endproclaim

\remark{Remarks}
The component operators $S_1,\dots,S_d$ of the $d$-shift 
are mutually commuting contractions in $\Cal B(H^2)$ which 
satisfy
$$
S_1S_1^* + \dots + S_dS_d^* = \bold 1 -E_0
$$
where $E_0$ is the projection onto the space of constant 
functions in $H^2$.  In particular, we conclude from 
Proposition 2.6 that for any 
$f_1,\dots,f_d\in H^2$
$$
\|S_1f_1+\dots + S_df_d\|^2 \leq \|f_1\|^2 + \dots + \|f_d\|^2.
$$

Notice too that if we replace $z_1,\dots,z_d$ with a different
set of coordinate functions $z_1^\prime,\dots,z_d^\prime$ 
for $\Bbb C^d$ then then the operators $(S_1,\dots,S_d)$ 
change to a new $d$-shift $(S_1^\prime,\dots,S_d^\prime)$.  
However, this change is not significant by virtue of the 
relation between $z_k$ and $z_k^\prime$.  More precisely, 
letting $V$ be the unitary operator defined on $\Bbb C^d$ 
by (2.5), one finds that 
$$
\Gamma(V)S_k\Gamma(V)^{-1} = S_k^\prime,\qquad k=1,2,\dots,d,   
$$
that is, $(S_1^\prime,\dots,S_d^\prime)$ {\it and} 
$(S_1,\dots,S_d)$ {\it are unitarily equivalent by way 
of a natural unitary automorphism of } $H^2$.  
In this sense we may speak of {\it the} $d$-shift acting 
on $H^2_d$.  In particular, we may conclude that each 
component operator $S_i$ is unitarily equivalent to every 
other one $S_j$, $1\leq j\leq d$.  

Finally, if $f$ is any polynomial in $\Cal P$ then we may
express $M_f$ as a polynomial in the operators $S_1,\dots,S_d$
as follows.  We find a polynomial 
function $g(w_1,\dots,w_d)$ of $d$ complex variables with the property
that $f$ is the composite function of $g$ with the coordinate 
functions $z_1,\dots,z_d$, 
$$
f(x) = g(z_1(x),\dots,z_d(x)),\qquad x\in B_d.  
$$
Once this is done the 
multiplication operator $M_f$ becomes the
corresponding polynomial in the operators $S_1,\dots,S_d$:
$$
M_f = g(S_1,\dots,S_d).  
$$
\endremark

We emphasize that in the higher dimensional cases $d\geq 2$, 
the operator norm $\|M_f\|$ can be larger than the sup norm
$\|f\|_\infty$ (see section 3 below).  
On the other hand, in all dimensions the spectral radius 
$r(M_f)$ of any polynomial multiplication operator satisfies
$$
r(M_f) = \sup_{z\in B_d}|f(z)|.  \tag{2.11}
$$ 
In the following result we establish the 
formula (2.11).  That follows from a straightforward 
application of the Gelfand theory of commutative Banach algebras
and we merely sketch the details.  

\proclaim{Proposition 2.12}
Let $\Cal A$ be the norm closed subalgebra of $\Cal B(H^2)$ generated
by the multiplication operators $M_f$, $f\in \Cal P$.  

Every element
of $\Cal A$ is a multiplication operator $M_f$ for some $f\in \Cal M$
which extends continuously to the closed ball $\bar B_d$,
and there is a natural homeomorphism of the closed unit ball onto 
the space $\sigma(\Cal A)$ of all complex homorphisms of $\Cal A$,  
$x\mapsto \omega_x$,  defined by 
$$
\omega_x(M_f) = f(x), \qquad \|x\|\leq 1.  
$$

For every such $f\in \Cal M$ one has 
$$
\lim_{n\to\infty}\|M_f^n\|^{1/n} = \sup_{\|x\|< 1}|f(x)|.  
$$
\endproclaim
\demo{proof}
Since the mapping $f\in \Cal M\mapsto M_f\in\Cal B(H^2)$ is an
isometric representation of the multiplier algebra on $H^2$ which
carries the unit of $\Cal M$ to that of $\Cal B(H^2)$, it is 
enough to work within $\Cal M$ itself.  That is, we may consider
$\Cal A$ to be the closure in $\Cal M$ of the algebra of polynomials,
and basically we need to identify its maximal ideal space.  

Because of the inequality $\|f\|_\infty \leq \|f\|_\Cal M$ of 
Proposition 2.2, we can assert that for every polynomial $f$ 
and every $x\in \Bbb C_d$ satisfying $\|x\|\leq 1$ we have 
$$
|f(x)| \leq \sup_{z\in B_d}|f(z)| =\|f\|_\infty\leq \|f\|_\Cal M.  
$$
It follows that there is a unique complex homomorphism
$\omega_x$ of $\Cal A$ satisfying 
$$
\omega_x(f) = f(x), \qquad f\in \Cal P.  
$$
For all $g\in \Cal A$ we now have a natural 
continuous extension $\tilde g$ 
of $g$ to the closed unit ball by setting  
$$
\tilde g(x) = \omega_x(g), \qquad \|x\|\leq 1.  
$$

$x\mapsto \omega_x$ is a one-to-one continuous
map of the closed ball in $\Bbb C^d$ 
onto its range in $\sigma(\Cal A)$.  
To see that it is surjective, let 
$\omega$ be an arbitrary element of $\sigma(\Cal A)$. 
Then for every $y\in \Bbb C^d$ we may consider the linear 
functional 
$$
\hat y(z) = \<z,y\>, \qquad z\in \Bbb C^d.  
$$
The map $y\mapsto \hat y$ is an antilinear mapping of $\Bbb C^d$ 
onto the space of linear functions in  $\Cal P$ 
and we claim that $\|\hat y\|_\Cal M\leq \|y\|$.  
Indeed, assuming that $y\neq 0$, 
the linear function 
$$
u(x) = \hat y(x)/\|y\| = \<x,y\>/\|y\|
$$
is part of a system of coordinates for $\Bbb C^d$.  
Proposition 2.6 implies $\|u\|_\Cal M \leq 1$, 
hence $\|\hat y\|_\Cal M \leq \|y\|$.  
Thus, $y\mapsto \omega(\hat y)$ defines an antilinear 
functional on $\Bbb C^d$ satisfying 
$$
|\omega(\hat y)| \leq \|\hat y\|_\Cal M \leq \|y\|,\qquad y\in \Bbb C^d.  
$$
It follows that there is a unique vector $x$ in the unit ball 
of $\Bbb C^d$ such that 
$$
\omega(\hat y) = \<x,y\>,\qquad y\in \Bbb C^d.  
$$
Thus, $\omega(f) = \omega_x(f)$ on every linear functional $f$.  
Since both $\omega$ and $\omega_x$ are continuous 
unital homomorphisms of $\Cal A$, since $\Cal P$ is the algebra
generated by the linear functions and the constants, 
and since $\Cal P$ is dense in 
$\Cal A$, it follows that $\omega = \omega_x$, and the 
claim is proved.  

Thus we have identified the maximal ideal space of $\Cal A$ with
the closed unit ball in $\Bbb C^d$.  From the elementary theory 
of commutative Banach algebras we deduce that for every 
$f$ in $\Cal A$, 
$$
\lim_{n\to \infty}\|f^n\|_\Cal M^{1/n} = r(f) = 
\sup\{|\omega(f)|: \omega\in\sigma(\Cal A)\} = 
\sup\{|\tilde f(x)|: \|x\|\leq 1\}=\|f\|_\infty,   
$$
completing the proof of Proposition 2.12\qed
\enddemo

The realization of the $d$-shift as a $d$-tuple of 
multiplication operators 
on the function space $H^2$ is not always convenient for 
making computations.  We require the following 
realization of $(S_1,\dots, S_d)$ as ``creation" operators on
the symmetric Fock space $\Cal F_+(E)$.  

\proclaim{Proposition 2.13}
Let $e_1,\dots e_d$ be an orthonormal basis for a 
Hilbert space $E$ of dimension $d$.  Define operators $A_1,\dots,A_d$ on 
$\Cal F_+(E)$ by 
$$
A_i\xi = e_i\xi, \qquad  \xi\in \Cal F_+(E)
$$
where $e_i\xi$ denotes the projection of $e_1\otimes \xi\in\Cal F(E)$ 
to the symmetric subspace $\Cal F_+(E)$.  Let $z_1,\dots,z_d$ be 
the system of orthogonal coordinates $z_i(x)=\<x,e_i\>$, 
$1\leq i\leq d$.  Then there is a unique unitary operator 
$W: H^2\to \Cal F_+(E)$ such that $W(1)=1$ and 
$$
W(z_{i_1}\dots z_{i_n}) = e_{i_1}\dots e_{i_n},\qquad 
n\geq 1, \quad i_k\in\{1,2,\dots,d\}. \tag{2.14}
$$
In particular, 
the $d$-tuple of operators $(A_1,\dots,A_d)$ is unitarily equivalent
to the $d$-shift.  
\endproclaim
\demo{proof}
For every $x\in E$ satisfying $\|x\|<1$ define an element 
$v_x\in \Cal F_+(E)$ by 
$$
v_x = 1\oplus x\oplus x^2\oplus x^3\oplus \dots.  
$$
it is obvious that $\|v_x\|^2 = (1-\|x\|^2)^{-1}$ and, more
generally, 
$$
\<v_x,v_y\> = (1-\<x,y\>)^{-1},\qquad \|x\|,\|y\|<1.  
$$
$\Cal F_+(E)$ is spanned by the set $\{v_x: \|x\|<1\}$.  

Let $\{u_x: \|x\|<1\}$ be the set of functions 
in $H^2$ defined in (1.11), and let $*$ be the unique 
conjugation of $E$ defined by $e_i^*=e_i$, that is, 
$$
(a_1e_1+\dots +a_de_d)^* = \bar a_1e_1+\dots+\bar a_de_d.  
$$
Then we have 
$$
\<u_x,u_y\> = (1-\<y,x\>)^{-1}=(1-\<x^*,y^*\>)^{-1} =
\<v_{x^*},v_{y^*}\>
$$
for all $x,y$ in the open unit ball of $E$.  By Proposition
1.12 there is a unique unitary operator $W:H^2\to \Cal F_+(E)$ 
such that $W(u_x)=v_{x^*}$, $\|x\|<1$.  

We have $W(1)=W(u_0) = v_0 = 1$.  Choose $x\in E$ 
satifying $\|x\|\leq 1$ and 
let $f_x$ denote the linear functional on $E$ defined by
$f_x(z)=\<z,x\>$.  We have $\|f_x\|_{H^2}\leq 1$ and in 
fact $\|f_x^n\|_{H^2}\leq 1$ for every $n=0,1,2,\dots$.  
Hence for every 
$0\leq r< 1 $ and every $z\in B_d$ we have 
$$
u_{r x}(z)=(1-\<z,r x\>)^{-1}=
\sum_{n=0}^\infty r^n\<z,x\>^n =\sum_{n=0}^\infty r^n f_x^n(z) \in H^2.  
$$
Similarly, 
$$
v_{r x^*} = \sum_{n=0}^\infty r^n (x^*)^n \in \Cal F_+(E).  
$$
Setting $W(u_{r x})$ equal to $v_{r x^*}$ and comparing coefficients 
of $r^n$ we obtain 
$$
W(f_x^n) = (x^*)^n
$$
for every $n=0,1,\dots$.  It follows that 
$$
W(f_{x_1}f_{x_2}\dots f_{x_n}) = x_1^*x_2^*\dots x_n^*
\tag{2.15}
$$
for every $x_1,x_2,\dots,x_n\in E$.  
Indeed, setting 
$$
\align
L(x_1,x_2,\dots,x_n) &= W(f_{x_1}f_{x_2}\dots f_{x_n})\\
R(x_1,x_2,\dots,x_n) &= x_1^*x_2^*\dots x_n^*
\endalign
$$
for $x_1,x_2,\dots,x_n\in E$, we see that both $L$ and $R$
are symmetric $n$-antilinear mappings which agree when
$x_1=x_2=\dots =x_n\in B_d$.  Hence $L=R$ and (2.15)
follows.  We obtain (2.14) by taking $x_k=e_{i_k}$ 
in (2.15).  

(2.14) obviously implies $WS_i = A_iW$ for $i=1,\dots,d$,
so that the $d$-tuples $(S_1,\dots,S_d)$ and $(A_1,\dots,A_d)$ 
are unitarily equivalent\qed
\enddemo


\subheading{3. von Neumann's inequality and the sup norm}

\proclaim{Definition 3.1}
A $d$-contraction is a $d$-tuple of operators 
$\bar T=(T_1,\dots,T_d)$ acting on a common Hilbert space
$H$ which commute with each other and satisfy 
$$
\|T_1\xi_1+\dots+T_d\xi_d\|^2 \leq \|\xi_1\|^2+\dots +\|\xi_d\|^2,
$$
for every $\xi_1,\dots,\xi_d\in H$.  
\endproclaim

\remark{Remark 3.2}
We make frequent use of the following observation.  
For operators 
$T_1,\dots,T_d$ on a common Hilbert space $H$,  
the following are equivalent.  
\roster
\item
$\|T_1\xi_1+\dots +T_d\xi_d\|^2 \leq \|\xi_1\|^2+\dots +\|\xi_d\|^2$
for all $\xi_1,\dots,\xi_d\in H$.  
\item
$T_1T_1^* + \dots +T_dT_d^* \leq \bold 1$.  
\endroster
To see this let $d\cdot H$ denote the direct sum of $d$ copies of 
$H$, and let $\bar T\in \Cal B(d\cdot H,H)$ be the operator 
defined by $\bar T(\xi_1,\dots,\xi_d) = T_1\xi_1+\dots+T_d\xi_d$.  
A simple computation shows that the adjoint 
$\bar T^*: H\to d\cdot H$ is given by
$$
\bar T^* \xi = (T_1^*\xi,\dots,T_d^*\xi).  
$$
Thus $\bar T\bar T^*$ is the operator in $\Cal B(H)$ given by
$\bar T\bar T^* = T_1T_1^*+\dots+T_dT_d^*$.  The equivalence
of (1) and (2) follows.  
\endremark

Notice that the $d$-shift $\bar S = (S_1,\dots,S_d)$ acting 
on $H^2_d$ is a $d$-contraction.  
Perhaps the most natural generalization of von Neumann's 
inequality for $d$-dimensional operator theory would 
make the following assertion.    
Let $\bar T = (T_1,\dots,T_d)$ be a $d$-contraction and let 
$f = f(z_1,\dots,z_d)$ be a polynomial in $d$ complex variables
$z_1,\dots,z_d$.  Then 
$$
\|f(T_1,\dots,T_d)\| \leq \sup_{\|z\|\leq 1}|f(z_1,\dots,z_d)|.  
$$
In this section we show that this inequality fails rather spectacularly 
for the $d$-shift, in that there is no constant $K$ for 
which 
$$
\|f(S_1,\dots,S_d)\| \leq K \sup_{\|z\|\leq 1}|f(z_1,\dots,z_d)|, 
$$
holds for all polynomials $f$.  It follows that the multiplier
algebra $\Cal M$ is a proper subalgebra of $H^\infty$.  Indeed, 
we exhibit continuous functions 
$$
f: \{z\in \Bbb C^d: \|z\|\leq 1\}\to \Bbb C
$$
which are analytic in the interior of the unit ball and which 
do not belong to $\Cal M$.  

We will establish the appropriate version 
of von Neumann's inequality for dimension $d \geq 2$
in section 8.  

\proclaim{Theorem 3.3} Assume $d\geq 2$.  
Let $c_0, c_1,\dots$ be a sequence of 
complex numbers having the properties 
$$
\align
&\sum_{n=0}^\infty |c_n| = 1, \tag{i}\\
&\sum_{n=0}^\infty |c_n|^2 n^{(d-1)/2} = \infty,  \tag{ii}
\endalign
$$
and define a function $f(z_1,\dots,z_d)$ for $|z_1|^2+\dots+|z_d|^2\leq 1$
as follows
$$
f(z_1,\dots,z_d) = \sum_{n=0}^\infty \frac{c_n}{s^n}(z_1z_2\dots z_d)^n, 
\tag{3.3}
$$
where $s$ denotes the sup norm 
$$
s = \sup_{|z_1|^2+\dots+|z_d|^2\leq 1} |z_1z_2\dots z_d| = \sqrt{\frac{1}{d^d}}.  
\tag{3.4}
$$
Then the power series (3.3) converges uniformly over the closed
unit ball to a function $f$ satisfying 
$\|f\|_\infty\leq 1$.  The restriction of $f$ to $B_d$ 
does not belong to $H^2$.  
Letting $f_0,f_1,f_2,\dots$ be the 
sequence of Taylor polynomials 
$$
f_N(z_1,\dots,z_d) = \sum_{n=0}^N \frac{c_n}{s^n}(z_1z_2\dots z_d)^n,  
$$
then we have $\|f_N\|_\infty\leq 1$ for every $N$ while
$$
\lim_{N\to \infty}\|f_N(S_1,\dots,S_d)\| =
\lim_{N\to\infty} \|f_N\|_\Cal M= \infty.  \tag{3.5}
$$
\endproclaim

\remark{Remarks}
It is clear that the function $f$ belongs to the ``ball algebra", 
that is the closure in the sup norm $\|\cdot\|_\infty$ of the algebra
of polynomials.  On the other hand, $f$ does not belong to the 
multiplier algebra $\Cal M$, and in particular the inclusion 
$\Cal M \subsetneq H^\infty$ is {\it proper}.  

Note too that it is a simple matter to give explicit examples of 
sequences $c_0,c_1,\dots$ satisfying conditions (i) and (ii).  For 
example, let $S$ be any {\it infinite} subset of the nonnegative 
integers which is sparse enough so that 
$$
\sum_{n\in S}\frac{1}{n^{(d-1)/4}}<\infty.  
$$ 
If we set $c_n = 1/n^{(d-1)/4}$ if $n\in S$ and $c_n=0$ otherwise,
then we obviously have (ii) because $S$ is infinite, and (i) can 
also be achieved after multiplying the sequence by a suitable 
positive constant.  
\endremark
\demo{proof}
The formula (3.4) for the sup norm, 
$$
s=d^{-d/2}
$$
follows from the elementary fact that
$$
(|z_1|^2|z_2|^2\dots |z_d|^2)^{1/d} \leq 
(|z_1|^2+|z_2|^2+\dots+|z_d|^2)/d
$$
with equality iff $|z_1|=|z_2|=\dots=|z_d|$. 

Let $p$ be the homogeneous polynomial 
$p(z_1,\dots,z_d)= z_1z_2\dots z_d$.  Then for every $n=0,1,2,\dots$ 
we have 
$$
\|p^n\|_\infty = \|p\|_\infty^n = d^{-nd/2} = s^n.  
$$
It follows that the power series (3.3)
$$
\sum_{n=0}^\infty \frac{c_n}{s^n}p(z_1,\dots,z_d)^n
$$
converges uniformly over the unit ball to $f$.  
Thus it remains to establish the condition (3.5).  

Now for any polynomial $g$ in the $d$ complex variables
$z_1,\dots,z_d$ we have 
$$
\|g(S_1,\dots,S_d)\| \geq \|g(S_1,\dots,S_d)1\|_{H^2} = \|g\|_{H^2}.  
$$
Thus it suffices to show that the sequence of 
Taylor polynomials $f_0, f_1,f_2,\dots $ defined by 
the partial sums of the series diverges in the $H^2$ norm, that is
$$
\sup_{N}\|f_N\|_{H^2}^2 =
\sum_{n=0}^\infty |c_n|^2\frac{\|p^n\|_{H^2}^2}{\|p^n\|_\infty^2}
= \infty.  \tag{3.6}
$$

In order to establish (3.6) we will show that there is a positive 
constant $A$ such that 
$$
\|p^n\|^2_{H^2}\geq Ad^{-nd}n^{(d-1)/2}, \tag{3.7}
$$
for all $n=1,2,\dots$.  In view of the fact that 
$\|p^n\|_\infty^2=d^{-nd}$ and the series $\sum_n |c_n|^2n^{(d-1)/2}$
diverges,  (3.6) will follow.  

The estimate (3.7) is based on the following computation.  Since 
the result is a statement about certain norms in the 
symmetric Fock space over $\Bbb C^d$, 
it is likely that the result of Lemma 3.8 can be 
found in the literature.  Since we are not aware of an appropriate 
reference and since the estimate (3.7) depends 
essentially on these formulas, 
we have provided the details.  

\proclaim{Lemma 3.8}Let $e_1,e_2,\dots,e_d$ be an orthonormal 
basis for $E=\Bbb C^d$.  Then for every 
$d$-tuple of nonnegative integers $k=(k_1,\dots,k_d)$ we have 
$$
\|e_1^{k_1}e_2^{k_2}\dots e_d^{k_d}\|_{E^{|k|}}^2 = 
\frac{k_1!k_2!\dots k_d!}{|k|!}
$$
where $|k| = k_1+k_2+\dots+k_d$.  
\endproclaim

\remark{Remark}
Regarding notation, we have written $e_1^{k_1}e_2^{k_2}\dots e_d^{k_d}$
for the projection of the vector 
$$
e_1^{\otimes k_1}\otimes e_2^{\otimes k_2}\otimes\dots
\otimes e_d^{\otimes k_d} \in E^{\otimes |k|}
$$
to the symmetric subspace $E^{|k|}\subseteq E^{\otimes |k|}$.  
\endremark
\demo{proof}
For $y_1,\dots,y_p\in E=\Bbb C^d$ we write $y_1y_2\dots y_p$ for 
the projection of 
$y_1\otimes y_2\otimes \dots\otimes y_p\in E^{\otimes p}$
to the symmetric subspace $E^p$.  Fixing $a\in E$, and $p\geq 1$ 
we have 
an associated ``creation operator" $A: E^{p-1}\to E^p$
defined by
$$
A(x_1x_2\dots x_{p-1}) =a x_1x_2\dots x_{p-1}, \qquad x\in E.  
$$
We claim first that for $p\geq 1$, the adjoint $A^*: E^p\to E^{p-1}$ 
is given by 
$$
A^*(y_1y_2\dots y_p) = 
\frac{1}{p}\sum_{k=1}^p \<y_k,a\>y_1\dots \hat y_k \dots y_p \tag{3.9}
$$
where $\hat y_k$ means that the term $y_k$ is missing from the 
symmetric tensor product.  
Indeed, if $\zeta$ denotes the right side of (3.9) then for every 
$x\in E$ we have 
$$
\align
\<\zeta,x^{p-1}\> &= 
\frac{1}{p}\sum_{k=1}^p\<y_1,x\>\dots 
\<y_{k-1},x\>\<y_k,a\>\<y_{k+1},x\>\dots \<y_p,x\> \\
&= \<y_1\otimes\dots\otimes y_p,p^{-1}(a\otimes x^{p-1}+
x\otimes a\otimes x^{p-2}\dots+x^{p-1}\otimes a)\>. 
\endalign 
$$
Since
$$
p^{-1}(a\otimes x^{p-1}+
x\otimes a\otimes x^{p-2}\dots+x^{p-1}\otimes a) = ax^{p-1}\in E^p,
$$
the right side of the preceding formula becomes 
$$
\<y_1\otimes\dots \otimes y_p,ax^{p-1}\> = \<y_1y_2\dots y_p,ax^{p-1}\> 
=\<A^*(y_1y_2\dots y_p),x^{p-1}\>.  
$$
(3.9) now follows because $E^{p-1}$ is spanned by vectors 
of the form $x^{p-1}$, $x\in E$.  

To prove Lemma 3.8 we proceed by induction on the total degree
$|k|$.  The formula is obvious for $|k|=0$.  Assuming that 
$|k|\geq 1$ and that the formula has been established for total 
degree $|k|-1$ then we may assume (after relabelling the basis 
vectors $e_1,\dots,e_d$ if necessary) that $k_1\geq 1$.  

Taking $a=e_1$ in (3.9) and noting that $\<e_1,e_1\>=1$
and $\<e_1,e_j\>=0$ if $j=2,\dots,d$, we find that 
$$
A^*(e_1^{k_1}e_2^{k_2}\dots e_d^{k_d}) =
\frac{k_1}{|k|}e_1^{k_1-1}e_2^{k_2}\dots e_d^{k_d},
$$
hence
$$
\align
\<e_1^{k_1}\dots e_d^{k_d},e_1^{k_1}\dots e_d^{k_d}\> &=
\<A^*(e_1^{k_1}\dots e_d^{k_d}),e_1^{k_1-1}e_2^{k_2}\dots e_d^{k_d}\>\\
&=\frac{k_1}{|k|}\|e_1^{k_1-1}e_2^{k_2}\dots e_d^{k_d}\|^2.  
\endalign
$$
The required formula now follows from the induction hypothesis\qed

\enddemo

Setting $k_1=k_2=\dots =k_d=n$ in Lemma 3.8, we obtain
$$
\|(e_1e_2\dots e_d)^n\|^2_{E^{nd}} = \frac{(n!)^d}{(nd)!}.  
$$
The right side is easily estimated using Stirling's formula
$$
n!\sim \sqrt{2\pi}\,n^{n+1/2}\,e^{-n}
$$
and after obvious cancellations we find that 
$$
\frac{(n!)^d}{(nd)!} \sim (\frac{2\pi^{d-1}}{d})^{1/2}d^{-nd}n^{(d-1)/2}.  
$$

In order to deduce (3.7) from the latter, choose an orthonormal 
basis $e_1,e_2,\dots,e_d$ for $\Bbb C^d$ so that 
$$
z_k(x) = \<x,e_k\>, \qquad k=1,2,\dots,d.  
$$
Then 
$$
\|(z_1z_2\dots z_d)^n\|_{H^2}^2 = \|(e_1e_2\dots e_d)^n\|^2_{E^{nd}}
$$
and (3.7) follows after choosing $A$ to be a positive number
appropriately smaller than $\sqrt{2\pi^{d-1}/d}$.  
That completes the proof of Theorem 3.3\qed
\enddemo

\remark{Remark 3.10}We recall that a  
$d$-tuple of commuting operators $\bar T=(T_1,T_2,\dots,T_d)$ on
a  Hilbert space $H$ is said to be {\it subnormal} if there is a 
commuting $d$-tuple of normal operators $\bar N=(N_1,N_2,\dots,N_d)$
on a larger Hilbert space $K\supseteq H$ such that 
$$
T_k = N_k\restriction_H,\qquad k=1,2,\dots,d.  
$$
The one-dimensional unilateral 
shift can be extended to a unitary operator on a larger
space.  That situation
is unique to dimension $1$, as we have 
\endremark

\proclaim{Corollary 1}
For every $d\geq 2$ the $d$-shift is not subnormal.  
\endproclaim
\demo{proof}
In Proposition 2.12 we identified the maximal ideal
space of the unital Banach algebra generated by the $d$-shift 
with the closed unit ball in $\Bbb C^d$.  In  particular,
for every polynomial $f$ 
the spectral radius of $f(S_1,\dots,S_d)$ is given by 
$$
r(f(S_1,\dots,S_d)) = \sup_{|z_1|^2+\dots+|z_d|^2\leq 1}|f(z_1,\dots,z_d)|.  
$$
If the $d$-shift were subnormal then $f(S_1,\dots,S_d)$ 
would be a subnormal operator for every polynomial $f$, hence
its norm would equal its spectral radius \cite{22, Problem 162},
contradicting Theorem 3.3.\qed\enddemo

The two most common Hilbert spaces associated with the unit ball $B_d$
arise from measures.  These are the spaces $H^2(\partial B_d)$ associated
with normalized surface measure on the boundary of $B_d$ and the space 
$H^2(B_d)$ associated with normalized volume measure on $B_d$ \cite{37}.  
It is reasonable to ask if the space $H^2$ can be associated with some 
measure on $\Bbb C^d$.  The answer is no because that would imply
that the $d$-shift is subnormal, contradicting Corollary 1.  The 
details are as follows.  

\proclaim{Corollary 2}
There is no positive measure $\mu$ on $\Bbb C^d$, 
$d\geq 2$, with the property 
that 
$$
\|f\|_{H^2}^2 = \int_{\Bbb C^d} |f(z)|^2\, d\mu(z) 
$$
for every polynomial $f$.  
\endproclaim
\demo{proof}
Suppose that such a measure $\mu$ did exist.  $\mu$ must be a 
probability measure because $\|1\|_{H^2}=1$, and it must have 
finite moments of all orders.  

We claim that $\mu$ must have compact support.  Indeed, if 
$f$ is any linear functional on $\Bbb C^d$ of the form 
$f(x) = \<x,e\>$ where $e$ is a unit vector of $\Bbb C^d$ 
then by Lemma 3.8 we have 
$$
\|f^n\|_{H^2}^2 = \|e^n\|^2_{E^n} = 1
$$
for every $n=1,2,\dots$.  Hence 
$$
\int_{\Bbb C^d}|f(z)|^{2n}\, d\mu(z) = 1.  
$$
Taking $2n$th roots we find that the function $f$ has 
norm $1$ when it is 
considered an element in the space 
$L^p(\Bbb C^d,\mu)$ for $p = 2,4,6,\dots$.  Letting $X$ 
be the closed support of the measure $\mu$ we find that 
$$
\sup_{z\in X}|f(z)| = 
\lim_{n\to\infty}(\int_{\Bbb C^d}|f|^{2n}\,d\mu)^{1/2n} = 1.  
$$
This proves that for every $z\in X$ and $e$ in the unit ball 
of $\Bbb C^d$ we have 
$$
|\<z,e\>|\leq 1, 
$$
thus $X$ must be contained in the closed unit ball
of $\Bbb C^d$.  

Now we simply view the $d$-shift as a $d$-tuple of 
multiplication operators in the space 
$L^2(\mu)$.  Here, $S_k$ is multiplication by 
$z_k$ acting on the closure (in $L^2(\mu)$) of the space 
of polynomials.  This $d$-tuple $(S_1,\dots,S_d)$ is obviously
subnormal, contradicting Corollary 1 above\qed
\enddemo

\remark{Remark 3.11}
In the conventional approach to dilation theory one seeks 
normal dilations for operators or sets of operators.  Theorem 
3.3 implies that this approach is inappropriate for $d$-contractions
and the unit ball of $\Bbb C^d$ in dimension greater than
one.  Indeed, if $(N_1,\dots,N_d)$ is a $d$-tuple of mutually 
commuting normal operators whose joint spectrum is contained 
in the closed unit ball of $\Bbb C^d$ then for every 
polynomial $f\in \Cal P$ we have 
$$
\|f(N_1,\dots,N_d)\|\leq \sup_{z\in B_d}|f(z)|.  
$$
Since Theorem 3.3 implies that 
there are polynomials $f$ for which the inequality 
$$
\|f(S_1,\dots,S_d)\|\leq \sup_{z\in B_d}|f(z)|
$$
fails, one cannot obtain such operators 
$f(S_1,\dots,S_d)$ by compressing $f(N_1,\dots,N_d)$ to 
any subspace.  Thus the $d$-shift cannot be dilated to a 
normal $d$-tuple having its spectrum in the closed unit ball.  
\endremark


\subheading{4. Maximality of the $H^2$ norm}

The purpose of this section is to show that 
in every dimension $d=1,2,\dots$ the $H^2$ 
norm is distinguished among all Hilbert seminorms defined
on the space $\Cal P$ of polynomials by being the 
{\it largest} Hilbert seminorm which is appropriate for 
operator theory on the unit ball of $\Bbb C^d$.  As a 
consequence, we show that the function space $H^2$ is 
contained in every other Hilbert space of analytic 
functions on the open unit ball which has these 
natural properties.  

\proclaim{Definition 4.1}Let $z_1,\dots,z_d$ be 
a system of coordinate functions on $\Bbb C^d$.  
A Hilbert seminorm $\|\cdot\|$ defined on the 
space $\Cal P$ of all polynomials 
is said to be contractive if for every $a\in \Bbb C$ and 
every $f_1,\dots,f_d\in \Cal P$ we have 
$$
\|a1 + z_1f_1+\dots z_df_d\|^2 \leq 
|a|^2 + \|f_1\|^2+\dots+\|f_d\|^2. 
$$
\endproclaim

\remark{Remarks}
Proposition 2.6 asserts that {\it the $H^2$ norm 
is a contractive norm on $\Cal P$}.  From Proposition
4.2 below it follows that the Hilbert norms defined 
on $\Cal P$ by both $H^2(B_d)$ and 
$H^2(\partial B_d)$ are contractive norms.  

It is a simple exercise to show that if a Hilbert
seminorm $\|\cdot\|$ is contractive relative to one system
of coordinates $z_1,\dots,z_d$ then it is contractive 
relative to every system of coordinates.  Thus the definition 
of contractive seminorm depends only on the structure of 
$\Bbb C^d$ as a $d$-dimensional Hilbert space.  
\endremark

Notice too that if $\|\cdot\|$ is any contractive seminorm
then for any system of coordinate functions $z_1,\dots,z_d$ 
the multiplication operators $(M_{z_1},\dots,M_{z_d})$ give 
rise to a $d$-contraction acting on the Hilbert space 
obtained by completing $\Cal P$ in this seminorm.  Indeed, we have 
the following somewhat more concrete characterization of 
contractive Hilbert seminorms.  

\proclaim{Proposition 4.2}
Let $\|\cdot\|$ be an arbitrary Hilbert seminorm on $\Cal P$, 
let $H$ be 
the inner product space defined by $\|\cdot\|$, 
and let $\Cal P_0$ be the maximal ideal in $\Cal P$ consisting 
of all polynomials $f$ such that $f(0)=0$.  Then $\|\cdot\|$ is 
a contractive seminorm iff the following two conditions are 
satisfied
\roster
\item
$1\perp\Cal P_0$ in the space $H$, and 
\item
For some system of coordinate functions $z_1,\dots,z_d$
the multiplication operators 
$(M_{z_1},\dots,M_{z_d})$ define a $d$-contraction on $H$.  
\endroster
\endproclaim
\demo{proof}
Once one notes that the most general element of $\Cal P_0$ is a sum of 
the form $z_1f_1+\dots +z_df_d$ with $f_1,\dots,f_d\in \Cal P$, 
the argument is straightforward\qed
\enddemo

We collect the following observation, which asserts that condition
(2) alone is enough in the presence of minimal symmetry.  

\proclaim{Corollary}
For every $\lambda$ in the circle group $\{z\in \Bbb C: |z|=1\}$
and every $f\in \Cal P$ set $f_\lambda(z)=f(\lambda z)$, 
$z\in \Bbb C^d$.  Let $\|\cdot\|$ be a Hilbert seminorm on 
$\Cal P$ which satisfies $\|f_\lambda\| = \|f\|$ for every 
$f\in \Cal P$ and every $\lambda$, 
such that for some system of coordinate functions 
$z_1,\dots,z_d$, the multiplication operators $M_{z_1},\dots,M_{z_d}$
give rise to a $d$-contraction acting 
on the Hilbert space $H$ obtained from $\|\cdot\|$.  

Then $\|\cdot\|$ is a contractive seminorm.  
\endproclaim
\demo{proof}
We show that the symmetry hypothesis implies
condition (1) of Proposition 4.2.  For every $\lambda$ in the unit 
circle we can define a unitary operator $U_\lambda$ uniquely on 
$H$ by setting 
$$
U_\lambda f = f_\lambda, \qquad f\in \Cal P.  
$$
It is obvious that $U$ is a unitary representation of the 
circle group on $H$.  Moreover, if $f$ is a homogeneous polynomial 
of degree $n=0,1,\dots$ then we have 
$$
U_\lambda f = \lambda^n f
$$
for all $\lambda$.  Thus for the inner product $\<\cdot,\cdot\>$ 
associated with $\|\cdot\|$ we have 
$$
\<f,1\> = \<U_\lambda f,U_\lambda 1\> = \lambda^n\<f,1\>, 
$$
so that if $n\geq 1$
then $\<f,1\>=0$.  It follows
that $1\perp \Cal P_0$, as required\qed
\enddemo

Following is the main result of this section.  

\proclaim{Theorem 4.3}
Let $\|\cdot\|$ be any contractive Hilbert seminorm on 
$\Cal P$.  Then for every $f\in \Cal P$ we have 
$$
\|f\| \leq k \|f\|_{H^2}
$$
where $k=\|1\|$.  In particular, the $H^2$ norm is the 
largest contractive Hilbert seminorm which assigns norm 
$1$ to the constant polynomial $f = 1$.  
\endproclaim

In particular, we see that the Hilbert norms arising from 
the ``Hardy" space $H^2(\partial B_d)$ and 
the ``Bergman" space $H^2(B_d)$ are both dominated by 
$\|\cdot\|_{H^2}$.  Indeed, we have the following inclusions
of the corresponding Hilbert spaces of analytic functions 
in the open ball $B_d$
$$
H^2\subseteq H^2(\partial B_d)\subseteq H^2(B_d),   
$$
where both inclusion maps are compact operators
of norm $1$.  Since 
we do not require the latter assertion, we omit the proof.  
However, note that of these three function spaces, $H^2$ 
is the only one that does not contain $H^\infty$, and it 
is the only one of the three for which the $d$-contraction 
defined by the multiplication operators 
$(M_{z_1},\dots,M_{z_d})$ {\it fails} to be subnormal.

\remark{Remark 4.4}
Every $d$-contraction $(T_1,\dots,T_d)$ in $\Cal B(H)$ 
gives rise to a normal completely positive map $P$ on $\Cal B(H)$ 
by way of 
$$
P(A) = T_1AT_1^* + \dots + T_dAT_d^*, \qquad A\in\Cal B(H).   
$$
Because of Remark 3.2 we have 
$P(\bold 1)=T_1T_1^*+\dots +T_dT_d^* \leq \bold 1$, 
and in fact the sequence
$A_n = P^n(\bold 1)$ is decreasing: 
$A_0=\bold 1\geq A_1\geq A_2\geq \dots \geq 0$.  Thus
$$
A_\infty = \lim_{n\to \infty}P^n(\bold 1)
$$
exists as a limit in the strong operator topology and 
satisfies $0\leq A_\infty \leq \bold 1$.  A $d$-contraction
$\bar T=(T_1,\dots,T_d)$ is called {\it null} if $A_\infty=0$.  
Notice that if the row norm of $\bar T$ is less than $1$, 
i.e., $T_1T_1^*+\dots+T_dT_d^*\leq r\bold 1$ for some $0<r<1$, 
then $\|P\|=\|P(\bold 1)\|\leq r<1$ and hence 
$\bar T$ is a null $d$-contraction.  
\endremark

For the proof of Theorem 4.3 we require an 
operator theoretic result which relates closely to the material
of section 6.  Gelu Popescu has pointed out that the operator 
$L$ is related to his Poisson kernel 
operator $K_r$ of \cite{35, section 8}, when $r=1$.

\proclaim{Theorem 4.5}
Let $(T_1,\dots,T_d)$ be a $d$-contraction on a Hilbert space
$H$, define the operator 
$$
\Delta = (\bold 1 - T_1T_1^* - \dots -T_dT_d^*)^{1/2}, 
$$
and the subspace $K=\overline{\Delta H}$.  Let $E$ be a 
$d$-dimensional Hilbert space and let 
$$
\Cal F_+(E) = \Bbb C\oplus E\oplus E^2\oplus\dots 
$$
be the symmetric Fock space over $E$.  

Then for every orthonormal basis $e_1,\dots, e_d$ for $E$ there 
is a unique bounded operator 
$L:\Cal F_+(E)\otimes K\to H$
satisfying $L(1\otimes \xi) = \Delta \xi$ and 
$$
L(e_{i_1}e_{i_2}\dots e_{i_n}\otimes \xi) = 
T_{i_1}T_{i_2}\dots T_{i_n}\Delta\xi
$$
for every $i_1,\dots,i_n\in\{1,2,\dots,d\}$, $n=1,2,\dots$.  
In general we have $\|L\|\leq 1$, 
and if $(T_1,\dots,T_d)$ is a null 
$d$-tuple, then $L$ is a coisometry: $LL^*=\bold 1_H$.  
\endproclaim

\demo{proof}
If there is a bounded operator $L$ satisfying the stated condition
then it is obviously unique because $\Cal F_+(E)$ is spanned by 
the set of vectors 
$$
\{1,e_{i_1},e_{i_2}e_{i_3}, e_{i_4}e_{i_5}e_{i_6},\dots:
i_k\in\{1,2,\dots,d\}, k=1,2,\dots\}.
$$  

We define $L$ by exhibiting its adjoint.  That is, we
will exhibit an operator $A: H\to \Cal F(E)\otimes K$, 
$\Cal F(E)$ denoting the full Fock space over $E$, and we will 
show that $\|A\|\leq 1$ and 
$$
\<L(\zeta),\eta\> = \<\zeta,A(\eta)\> \tag{4.6}
$$
for $\zeta$ of the form $1\otimes \xi$ or 
$e_{i_1}e_{i_2}\dots e _{i_n}\otimes \xi$ for $\xi\in K$.  
At that 
point we can {\it define} $L$ to be the adjoint of $P_+A$, 
$P_+$ denoting the projection of $\Cal F(E)$ onto its 
subspace $\Cal F_+(E)$.  

For every $\eta\in H$, we define $A\eta$ as a sequence of 
vectors $(\zeta_0,\zeta_1,\zeta_2,\dots)$ where 
$\zeta_n\in E^{\otimes n}\otimes K$ is defined by 
$$
\zeta_n=\sum_{i_1,\dots,i_n=1}^d e_{i_1}\otimes \dots\otimes e_{i_n}
\otimes \Delta T_{i_n}^*\dots T_{i_1}^*\eta 
$$
for $n\geq 1$ and $\zeta_0 = 1\otimes \Delta\eta$.  
Notice that since $T_1^*,\dots,T_d^*$ commute, $\zeta_n$
actually belongs to the symmetric subspace 
$E^n\otimes K$.  We claim first that 
$$
\sum_{n=0}^\infty \|\zeta_n\|^2 \leq \|\eta\|^2
$$
so that in fact $A$ maps into $\Cal F(E)\otimes K$ and is 
a contraction.  Indeed, we have 
$$
\|\zeta_n\|^2 = 
\sum_{i_1,\dots,i_n=1}^d\|\Delta T_{i_n}^*\dots T_{i_1}^*\eta\|^2
= \sum_{i_1,\dots,i_n=1}^d\<T_{i_1}\dots T_{i_n}
\Delta^2T_{i_n}^*\dots T_{i_1}^*\eta,\eta\>.  
$$
Let $P(A) = T_1AT_1^*+\dots+T_dAT_d^*$ 
be the completely positive map of Remark 4.4.  Noting 
that $\Delta^2 = \bold 1-P(\bold 1)$ we find that 
$$
\sum_{i_1,\dots,i_n=1}^dT_{i_1}\dots 
T_{i_n}\Delta^2T_{i_n}^*\dots T_{i_1}^* = P^n(\bold 1-P(\bold 1)) =
P^n(\bold 1) - P^{n+1}(\bold 1),
$$
and hence
$$
\|\zeta_n\|^2 = \<P^n(\bold 1)\eta,\eta\> - \<P^{n+1}(\bold 1)\eta,\eta\>.  
$$
The series $\|\zeta_0\|^2+\|\zeta_1\|^2+\dots $therefore 
telescopes and we are left with  
$$
\sum_{n=0}^\infty \|\zeta_n\|^2 = \|\eta\|^2 - \<A_\infty \eta,\eta\>
\leq \|\eta\|^2, \tag{4.7}
$$
where $A_\infty$ is the positive contraction 
$A_\infty = \lim_{n\to \infty }P^n(\bold 1)$
of Remark 4.4.  

We now verify (4.6) for $\zeta$ of the form 
$\zeta = e_{j_1}\dots e_{j_n}\otimes \xi$ for $n\geq 1$, 
$j_1,\dots,j_n\in \{1,2,\dots,d\}$ and 
$\xi\in K$.  We have 
$$
\align
\<e_{j_1}\dots e_{j_n}\otimes \xi,A\eta\>
&= \sum_{i_1,\dots,i_n=1}^d\<e_{j_1}\dots e_{j_n}\otimes \xi,
e_{i_1}\otimes\dots \otimes e_{i_n}\otimes 
\Delta T_{i_n}^*\dots T_{i_1}^*\eta\> \\
&= \sum_{i_1,\dots,i_n=1}^d \<e_{j_1}\otimes\dots\otimes e_{j_n},
e_{i_1}\otimes\dots\otimes e_{i_n}\>
\<\xi,\Delta T_{i_n}^*\dots T_{i_1}^*\eta\> \\
&= \<\xi,\Delta T_{j_1}^*\dots T_{j_1}^*\eta\> = 
\<T_{j_1}\dots T_{j_n}\Delta\xi,\eta\> = \< L(\zeta),\eta\>.  
\endalign
$$
For $\zeta=1\otimes \xi$ with $\xi\in K$ we have 
$$
\<1\otimes \xi,A\eta\> = \<1\otimes \xi,1\otimes\Delta\eta\> 
=\<\xi,\Delta \eta\> = \<\Delta\xi,\eta\>,
$$
as required.  If $(T_1,\dots,T_d)$ is a null 
$d$-tuple, then it is clear from (4.7) that $A$ is 
an isometry and hence $L$ is a coisometry.  
\qed
\enddemo

\demo{proof of Theorem 4.3}
Let $H$ be the Hilbert space obtained by completing $\Cal P$ 
in the seminorm $\|\cdot\|$.  Choose an orthonormal basis 
$e_1,\dots,e_d$ for $E=\Bbb C^d$ 
and let $z_1,\dots,z_d$ be the corresponding system 
of coordinate functions $z_i(x) = \<x,e_i\>$, $i=1,\dots,d$.  

Since $\|\cdot\|$ is a contractive Hilbert seminorm the 
multiplication operators 
$$
T_k = M_{z_k},\qquad k=1,\dots,d
$$
define a $d$-contraction $(T_1,\dots,T_d)$ in $\Cal B(H)$.  
Set 
$$
\Delta = (\bold 1-\sum_{k=1}^d T_kT_k^*)^{1/2}, 
$$
let $K=\overline{\Delta H}$ be the closed range of $\Delta$ 
and let $L: \Cal F_+(E)\otimes K\to H$ be the contraction 
defined in Theorem 4.5 by the conditions 
$L(1\otimes \xi)=\Delta \xi$ and, for $n=1,2,\dots$ 
$$
L(e_{i_1}\dots e_{i_n}\otimes \xi) = T_{i_1}\dots T_{i_n}\Delta \xi,
\tag{4.8}
$$
$\xi\in K$, $i_1,\dots,i_n\in \{1,2,\dots,d\}$.  

The constant polynomial $1\in \Cal P$ is represented by a vector 
$v$ in $H$.  We claim that $\Delta v=v$.  Indeed, since 
$\|\cdot\|$ is a contractive seminorm, condition (1) of 
Proposition 4.2 implies that 
$$
v\perp T_1H + T_2H+\dots +T_dH,
$$
and hence $T_k^*v=0$ for $k=1,\dots,d$.  It follows that 
$$
\|\Delta v\|^2 =\<\Delta^2v,v\> = \|v\|^2 -\sum_{k=1}^d\|T_k^*v\|^2 
= \|v\|^2,
$$
and hence $\Delta v=v$ because $0\leq \Delta\leq \bold 1$.  

In particular, $v=\Delta v\in \overline{\Delta H} = K$.  Taking
$\xi = v$ in (4.8) we obtain
$$
L(e_{i_1}\dots e_{i_n}\otimes v) = T_{i_1}\dots T_{i_n}v.  
$$
Since $v$ is the representative of $1$ in $H$, 
$T_{i_1}\dots T_{i_n}v$ is the representative of the 
polynomial $z_{i_1}\dots z_{i_n}$ in $H$, and we have 
$$
L(e_{i_1}\dots e_{i_n}\otimes v) = z_{i_1}\dots z_{i_n}\in H.  
$$

By Propostion 2.13 there is a unitary operator 
$W:H^2\to \Cal F_+(E)$ which carries $1$ to $1$ and carries 
$z_{i_1}\dots z_{i_n}\in H^2$ to 
$e_{i_1}\dots e_{i_n}\in \Cal F_+(E)$. Hence 
$$
L(W(z_{i_1}\dots z_{i_n})\otimes v)=z_{i_1}\dots z_{i_n}.  
$$
By taking linear combinations we find that for every 
polynomial $f\in \Cal P$ 
$$
L(Wf\otimes v) = f
$$
where, $f$ on the left is considered an element of $H^2$ 
and $f$ on the right is considered an element of $H$.  Since
$\|L\|\leq 1$ and $W$ is unitary, we immediately 
deduce that 
$$
\|f\|_H \leq \|Wf\otimes v\| = \|f\|_{H^2}\cdot \|v\|_H.  
$$
Theorem 4.3 follows after noting that $\|v\|_H=\|1\|_H$\qed
\enddemo

\remark{Remarks}
In particular, the $H^2$ norm is the {\it largest} Hilbert 
seminorm $\|\cdot\|$ on the space $\Cal P$ of all polynomials
which is contractive and is normalized
so that $\|1\|= 1$.  
\endremark

The following result gives a precise sense 
in which the space $H^2$ will 
be present in any function space that is 
appropriate for doing operator theory relative to the 
unit ball in $\Bbb C^d$.

\proclaim{Corollary}
Let $\Cal F$ be a linear space of analytic functions defined 
on the open unit ball in $\Bbb C^d$, $d=1,2\dots$, which 
contains all polynomials and is a Hilbert 
space relative to some Hilbert norm with the following two 
properties:
\roster
\item
The constant functions are orthogonal to the space $\Cal P_0$ of 
polynomials which vanish at $0$
\item
For some orthonormal set of linear functionals 
$f_1,\dots,f_d$ on $\Bbb C^d$ the multiplication operators 
$(M_{f_1},\dots,M_{f_d})$ define a $d$-contraction on this 
Hilbert space.  
\endroster

Then $\Cal F$ contains all functions in $H^2$ and the inclusion
map of $H^2$ into $\Cal F$ is a bounded operator.  
\endproclaim

\remark{Remark}
By an orthonormal set of linear functionals we of course mean
a set of linear functionals of the form $f_k(x)= \<x,e_k\>$ 
where $e_1,\dots,e_d$ is an orthonormal basis for $\Bbb C^d$.  
\endremark
\demo{proof}
Let $H$ be the subspace of $\Cal F$ obtained by closing 
the space of polynomials in the norm $\|\cdot\|$ of $\Cal F$. 
Since $f\in \Cal P\mapsto \|f\|$ is a contractive Hilbert norm 
on polynomials, Theorem 4.3 implies that $\|f\|\leq k \|f\|_{H^2}$
for every $f\in \Cal P$, where $k=\|1\|>0$.  Thus $H^2$ is 
contained in $H\subseteq \Cal F$ and the inclusion of 
$H^2$ in $\Cal F$ is of norm at most $k$\qed
\enddemo

We will make use of the following extremal property
of the $H^2$ norm below.  

\proclaim{Corollary 4.9}
Let $\|\cdot\|$ be a contractive Hilbert seminorm on 
$\Cal P$ satisfying $\|1\|=1$ 
and let $z_1,\dots,z_d$ be a system of 
orthogonal coordinate functions for $E=\Bbb C^d$.  
Then for every $n=1,2,\dots$ we have 
$$
\sum_{i_1,\dots,i_n=1}^d \|z_{i_1}z_{i_2}\dots z_{i_n}\|^2 
\leq \frac{(n+d-1)!}{n!(d-1)!}, \tag{4.10}
$$
with equality holding iff $\|f\| = \|f\|_{H^2}$ for 
every polynomial $f$ of degree at most $n$.  
\endproclaim

\demo{proof}
Let $\{e_1,\dots,e_d\}$ be an orthonormal basis 
for a $d$-dimensional Hilbert space $E$.  We consider the 
projection $P_n\in\Cal B(E^{\otimes n})$ of the 
full tensor product onto its symmetric subspace $E^n$.  
Since $\|\cdot\|$ is a contractive seminorm, Theorem 
4.3 implies that for all $i_1,\dots,i_n$ we have 
$$
\|z_{i_1}\dots z_{i_n}\|\leq \|z_{i_1}\dots z_{i_n}\|_{H^2}
= \|P_n(e_{i_1}\otimes\dots\otimes e_{i_n})\|
$$
and hence
$$
\sum_{i_1,\dots,i_n=1}^d\|z_{i_1}\dots z_{i_n}\|^2 \leq 
\sum_{i_1,\dots,i_n=1}^d \|P_n(e_{i_1}\otimes\dots\otimes e_{i_n})\|^2.  
$$
Since 
$\{e_{i_1}\otimes\dots\otimes e_{i_n}: 1\leq i_1,\dots,i_n\leq d\}$
is an orthonormal basis for $E^{\otimes n}$ the term on the 
right is $\tr(P_n) = \dim(E^n)$, and (4.10) follows from 
the computation of the dimension of $E^n$ in A.5.  

Let $\Cal P_n$ denote the subspace of $H^2$ consisting of
homogeneous polynomials of degree $n$, and let $Q_n$ be the 
projection of $H^2$ on $\Cal P_n$.  The preceding 
observations imply that if $A$ is any operator on 
$H^2$ which is supported in $\Cal P_n$ in the sense that 
$A=Q_nAQ_n$ then the trace of $A$ is given by
$$
\tr(A) = \sum_{i_1,\dots,i_n=1}^d
\<Az_{i_1}\dots z_{i_n},z_{i_1}\dots z_{i_n}\>_{H^2}.  \tag{4.11}
$$
Now fix $n$ and suppose equality holds in (4.10).  Since $\|\cdot\|$
is a contractive Hilbert seminorm satisfying $\|1\|=1$, Theorem 
4.3 implies that $\|f\|\leq \|f\|_{H^2}$ for every 
$f\in\Cal P$, and hence there is a unique operator $H\in\Cal B(H^2)$
satisfying 
$$
\<f,g\> = \<Hf,g\>_{H^2}, \qquad f,g\in\Cal P
$$
and one has $0\leq H\leq \bold 1$.  Considering the 
compression $Q_nHQ_n$ of $H$ to $\Cal P_n$ we see from 
(4.11) that 
$$
\tr(Q_nHQ_n)=\dim(E^n)=\tr(Q_n).  
$$
Since $Q_n-Q_nHQ_n\geq 0$ and the trace is faithful, we 
conclude that $Q_nHQ_n=Q_n$, and since $H$ is a positive 
contraction it follows that $Hf=f$ for every $f\in\Cal P_n$.  

We claim that $Hf=f$ for every $f\in\Cal P_k$ and every 
$k=0,1,\dots,n$.  To see that, choose a linear functional
$z\in\Cal P$ satisfying $\|z\|_{H^2}=1$.  Since $\|\cdot\|$
is a contractive seminorm we have $\|z\cdot f\|\leq \|f\|$ 
for every $f\in\Cal P$ and in particular we have 
$\|z^n\|=\|z^{n-k}z^k\|\leq \|z^k\|$.  Thus
$$
\<Hz^k,z^k\>_{H^2}=\|z^k\|^2 \geq \|z^n\|^2=\|z^n\|^2_{H^2}.  
$$
Since the $H^2$ norm of any power of $z$ is $1$ and 
 $0\leq H\leq \bold 1$, 
it follows that $H z^k=z^k$.  
Since every polynomial of degree at most $n$ is a linear 
combination of monomials of the form $z^k$ with $z$ 
as above and $k=0,1,\dots,n$, the proof of Theorem 4.9
is complete \qed
\enddemo


\subheading{5.  The Toeplitz $C^*$-algebra}
Let $\bar S = (S_1,\dots,S_d)$ be the $d$-shift.  

\proclaim{Definition 5.1}
The Toeplitz \cstar\ is the \cstar\ $\Cal T_d$ generated by 
the operators $S_1,\dots,S_d$.  
\endproclaim

\remark{Remarks}
Notice that we have not included the identity operator 
as one of the generators of $\Cal T_d$, so that $\Cal T_d$ 
is by definition the norm closed linear span of the set 
of finite products of the form $T_1T_2\dots T_n$, 
$n=1,2,\dots$, where 
$$
T_i\in \{S_1,\dots,S_d,S_1^*,\dots, S_d^*\}.
$$
Nevertheless, (5.5) below implies that $\Cal T_d$ contains an
invertible positive operator
$$
(d\bold 1+N)(\bold 1+N)^{-1} = S_1^*S_1+\dots+S_d^*S_d,
$$  
and hence $\bold 1\in \Cal T_d$.  Thus $\Cal T_d$ is 
the \cstar\ generated by all multiplication operators 
$M_f\in \Cal B(H^2)$, $f\in \Cal P$.  
\endremark

If one starts with the Hilbert space $H^2(\partial B_d)$ rather
than $H^2$ then there is a natural Toeplitz \cstar 
$$
\Cal T_{\partial B_d} = C^*\{M_f: f\in \Cal P\}\subseteq 
\Cal B(H^2(\partial B_d)),
$$
and similarly there is a Toeplitz \cstar\ $\Cal T_{B_d}$ on 
the Bergman space 
$$
\Cal T_{B_d} = C^*\{M_f: f\in \Cal P\}\subseteq \Cal B(H^2(B_d))
$$
\cite{16}.  In fact, it is not hard to 
show that the three \cstar s $\Cal T_d$, 
$\Cal T_{\partial B_d}$ and $\Cal T_{B_d}$ are unitarily 
equivalent.  In that sense, the \cstar\ $\Cal T_d$ is not new.  
However, we are concerned with the relationship between 
the $d$-shift and its enveloping \cstar\ $\Cal T_d$, 
and here there are some
essential differences.  

For example, in the classical case of $H^2(\partial B_d)$ one 
can start with a continuous complex-valued function 
$f\in C(\partial B_d)$ and define a Toeplitz operator $T_f$ 
on $H^2(\partial B_d)$ by compressing the operator of multiplication 
by $f$ (acting on $L^2(\partial B_d)$) to the subspace 
$H^2(\partial B_d)$.  In our case however, continuous symbols 
do not give rise to Toeplitz operators.  Indeed, we have seen 
that there are continuous functions $f$ on the closed unit ball 
which are uniform limits of holomorphic polynomials, 
but which do not belong 
to $H^2$.  For such an $f$ the ``Toeplitz" operator $T_f$ is 
not defined.  
Thus we have taken some care to develop the properties 
of $\Cal T_d$ that we require.  

Let $N$ be the number operator acting on $H^2$, defined as 
the generator of the one-parameter unitary group 
$$
\Gamma(e^{it}\bold 1_E) = e^{itN}, \qquad t\in \Bbb R, 
$$
$\Gamma$ being the representation of the unitary group 
of $E$ on $H^2$ defined in the remarks
following Definition 2.10.  
$N$  obviously has discrete spectrum $\{0,1,2,\dots\}$ and 
the $n$th eigenspace of $N$ is the space $\Cal P_n$ of 
homogeneous polynomials of degree $n$,
$$
\Cal P_n = \{\xi\in H^2: N\xi = n\xi\}, \qquad n=0,1,2,\dots.  
$$
$(\bold 1+N)^{-1}$ is a compact operator, 
and it is a fact that for every real number $p>0$,
$$
\tr (\bold 1+N)^{-p} <\infty \iff p>d.  \tag{5.2}
$$
Since $N$ is unitarily equivalent to the Bosonic number operator, 
the assertion (5.2) is probably known.  We lack an appropriate 
reference, however, and have included a proof of (5.2) in 
the Appendix for the reader's convenience.  

The following result exhibits the commutation relations 
satisfied by the $d$-shift. 
 
\proclaim{Proposition 5.3}
Suppose that $d=2,3,\dots$ and let $(S_1,\dots,S_d)$ be 
the $d$-shift.  Then for all $i,j=1,\dots,d$ we have 
$$
S_i^*S_j-S_jS_i^* = (\bold 1 + N)^{-1}(\delta_{ij}\bold 1-S_jS_i^*),
\tag{5.4}
$$
and
$$
S_1^*S_1+\dots+S_d^*S_d = (d\bold 1+N)(\bold 1+N)^{-1}.  
\tag{5.5}
$$
In particular, $\|S_1^*S_1+\dots +S_d^*S_d\|=d$.  The 
commutators $S_i^*S_j-S_jS_i^*$ belong to every Schatten 
class $\Cal L^p(H^2)$ for $p>d$, but they 
do not belong to $\Cal L^d(H^2)$.  
\endproclaim

\remark{Remark}
It follows that if $A,B$ are operators belonging to the 
unital $*$-algebra generated by $S_1,\dots,S_d$, then 
$AB-BA\in \Cal L^p(H^2)$ for every $p>d$, and hence 
any product of at least $d+1$ such commutators belongs
to the trace class.  
\endremark

\demo{proof}To establish these formulas 
it is more convenient to work with the $d$-shift in its 
realization on $\Cal F_+(E)$ described in Proposition 2.13.  
Thus, we pick an orthonormal basis $e_1,\dots,e_d$ for a 
$d$-dimensional Hilbert space $E$ and set 
$$
S_i\xi = e_i\xi, \qquad 1\leq i\leq d, 
$$
for $\xi\in \Cal F_+(E)= \Bbb C\oplus E\oplus E^2\oplus\dots$.
The number operator $N$ acts as follows on $E^n$: 
$$
N\xi = n\xi, \qquad \xi\in E^n, n=0,1,2,\dots.  
$$  

We first establish (5.4).  It suffices to verify that the operators 
on both sides of (5.4) agree on every finite dimensional space 
$E^n$, $n=0,1,2,\dots$.  For $n=0$ and $\lambda\in \Bbb C$ we 
have $S_i^*S_j\lambda=\lambda S_i^*e_j = \delta_{ij}\lambda$, 
while $S_jS_i^*\lambda = 0$, hence (5.4) holds on $\Bbb C$.  
For $n\geq 1$ and $\xi\in E^n$ of the form $\xi = y^n$ we see 
from formula (3.9) that 
$$
S_j^*S_i\xi = S_j^*(e_iy^n) = \frac{\delta_{ij}}{n+1}y^n +
\frac{n}{n+1}\<y,e_j\>e_iy^{n-1}
$$
while 
$$
S_iS_j^*\xi = \<y,e_j\>S_iy^{n-1}=\<y,e_j\>e_iy^{n-1}.  
$$
Hence 
$$
S_j^*S_i\xi -S_iS_j^*\xi= \frac{1}{n+1}(\delta_{ij}\xi-S_jS_i^*\xi).
$$
The latter holds for all $\xi\in E^n$ because $E^n$ is spanned 
by $\{y^n: y\in E\}$, and (5.4) follows.  

Formula (5.5) follows from (5.4).  Indeed, for 
$\xi\in E^n$ we have 
$$
S_i^*S_i\xi = \frac{1}{n+1}\xi + \frac{n}{n+1}S_iS_i^*\xi.  
$$
By the remarks following (2.10) we have 
$$
S_1S_1^* + \dots + S_dS_d^* = \bold 1 - E_0,   \tag{5.6}
$$
$E_0$ denoting the projection on $\Bbb C$.  Summing 
the previous formula on $i$ we obtain
$$
\sum_{i=1}^d S_i^*S_i\xi = \frac{d}{n+1}\xi + \frac{n}{n+1}(\xi-E_0\xi) 
= \frac{n+d}{n+1}\xi - \frac{n}{n+1}E_0\xi = \frac{n+d}{n+1}\xi, 
$$
and (5.5) follows.  

Now suppose $p>d$.  Because of (5.2) the operator
$(\bold 1+N)^{-1}$ belongs to $\Cal L^p$; since
$\Cal L^p$ is an ideal 
(5.4) implies that 
$S_i^*S_j-S_jS_i^*\in \Cal L^p$ for all $i,j$.  

Finally, we claim that no self-commutator 
$[S_i^*,S_i]=S_i^*S_i-S_iS_i^*$ belongs to $\Cal L^d$.  Indeed, 
since the operators $S_1,\dots,S_d$ are unitarily equivalent
to each other (by the remarks following Definition 2.10), we 
see that if one $[S_i^*,S_i]$ belongs to $\Cal L^d$ then they 
all do, and in that case we would have 
$$
\sum_{i=1}^d[S_i^*,S_i] \in \Cal L^d.  
$$
By (5.5) and (5.6) the left side of this formula is 
$$
\align
\sum_{i=1}^d S_i^*S_i - \sum_{i=1}^d S_iS_i^* &=
(d\bold 1 + N)(\bold 1+N)^{-1} - (\bold 1-E_0) \\
&= E_0 + (d-1)(\bold 1+N)^{-1}.  
\endalign
$$
Since $(\bold 1+N)^{-1}\notin \Cal L^d$ by (5.2), we 
have a contradiction and the proof of Proposition 5.3
is complete\qed
\enddemo

\remark{The $d$-shift and the Canonical Commutation Relations}
The $d$-shift is closely related to the creation operators 
$(C_1,\dots,C_d)$ associated with the canonical commutation 
relations for $d$ degrees of freedom.  Indeed, one can think
of $\bar S=(S_1,\dots,S_d)$ as the partial isometry 
occurring in the polar decomposition of $\bar C=(C_1,\dots,C_d)$
in the following way.  Choose an orthonormal basis 
$e_1,\dots,e_d$ for a $d$-dimensional Hilbert space $E$.  
For $k=1,\dots,d$, $C_k$ is defined 
on the dense subspace of $\Cal F_+(E)$ spanned by 
$E^n$, $n=0,1,\dots$ as follows 
$$
C_k\xi=\sqrt{n+1}e_k\xi,\qquad \xi\in E^n
$$
(see \cite{40}).  The $C_k$ are of course unbounded 
operators, and they satisfy the complex form of the 
canonical commutation relations
$$
C_i^*C_j-C_jC_i^* = \delta_{ij}\bold 1,\qquad 1\leq i,j\leq d.  
$$
One finds that the row operator
$$
\bar C=(C_1,\dots,C_d):
\underbrace{\Cal F_+(E)\oplus\dots\oplus\Cal F_+(E)}
_{\text{$d$ times}}\to\Cal F_+(E)
$$
is related to the number operator $N$ by 
$\bar C \bar C^*=N$, and in fact the polar decomposition 
of $\bar C$ takes the form 
$$
\bar C = N^{1/2}\bar S
$$
where $\bar S=(S_1,\dots,S_d)$ is the $d$-shift; i.e., 
$C_k=N^{1/2}S_k$, $k=1,\dots,d$.  
\endremark

We have seen that the $d$-shift is not a subnormal 
$d$-tuple.  The following result asserts that, at 
least, the individual operators $S_k$, $k=1,\dots,d$ 
are hyponormal.  Indeed, any linear combination of 
$S_1,\dots,S_d$ is a hyponormal operator.  

\proclaim{Corollary}
For every $k=1,\dots,d$ we have $S_k^*S_k\geq S_kS_k^*$.  
\endproclaim
\demo{proof}
Proposition 5.3 implies that 
$$
S_k^*S_k-S_kS_k^* = (\bold 1+N)^{-1}(\bold 1 - S_kS_k^*).  
$$
Since $\|S_k\|\leq 1$, both factors on the right are positive 
operators.  Let $E_n$ be the $n$th spectral projection 
of $N$, $n=0,1,\dots$.  Since $S_kE_n = E_{n+1}S_k$ 
it follows that $S_kS_k^*$ commutes with $E_n$.  
Thus $(\bold 1+N)^{-1}$ commutes with 
$\bold 1 - S_kS_k^*$, and the assertion follows.  
\qed
\enddemo

Of course in dimension $d=1$, the commutator 
$S^*S-SS^*$ is a rank-one operator and therefore 
belongs to every Schatten class $\Cal L^p$, $p\geq 1$.  

\proclaim{Theorem 5.7}
$\Cal T_d$ contains the algebra $\Cal K$ of all compact 
operators on $H^2$, and we have an exact sequence of 
\cstar s 
$$
0 @>>> \Cal K \hookrightarrow \Cal T_d @>>\pi> C(\partial B_d)@>>>0
$$
where $\pi$ is the unital $*$-homomorphism defined by 
$$
\pi(S_k) = z_k, 
$$
$z_k$ being the $k$th coordinate function $z_k(x)=\<x,e_k\>$, 
$x\in\partial B_d$.  

Letting $\Cal A$ be the commutative algebra of polynomials in the 
operators $S_1,\dots, S_d$ we have 
$$
\Cal T_d = \overline{\text{\rm span}} \Cal A\Cal A^*.  \tag{5.8}
$$
\endproclaim
\demo{proof}
Let $E_0$ be the one-dimensional projection onto the space 
of constants in $H^2$.  By the remark following Definition 
2.10 we have 
$$
E_0 = \bold 1- S_1S_1^* -\dots -S_dS_d^*\in
{\text{\rm span}}\Cal A\Cal A^*.  
$$
Thus for any polynomials $f,g$, the the rank-one 
operator
$$
f\otimes \bar g: \xi \mapsto \<\xi,g\>f
$$
can be expressed as 
$$
f\otimes \bar g = M_fE_0M_g^* \in {\text{\rm span}}\Cal A\Cal A^*.
$$
It follows that the norm closure of ${\text{\rm span}}\Cal A\Cal A^*$
contains the algebra $\Cal K$ of all compact operators.  

By Proposition 5.3, the quotient $\Cal T_d/\Cal K$ is a commutative
\cstar\ which is generated by commuting normal elements $Z_k=\pi(S_k)$,
$k=1,\dots,d$ satisfying 
$$
Z_1Z_1^*+\dots+Z_dZ_d^* = \bold 1.  
$$  
Because $\Cal T_d$ is commutative modulo $\Cal K$ and since
$\overline{\text{\rm span}}\Cal A\Cal A^*$ contains $\Cal K$, 
it follows that $\overline{\text{\rm span}}\Cal A\Cal A^*$ is 
closed under multiplication and (5.8) follows.

Let $X$ be the joint spectrum of the commutative normal 
$d$-tuple $(Z_1,\dots,Z_d)$ that generates $\Cal T_d/\Cal K$.  
$X$ is a nonvoid subset of the sphere $\partial B_d$, and 
we claim that $X=\partial B_d$.  Indeed, 
since the unitary group $\Cal U(E)$ acts transitively on 
$\partial B_d$ it suffices to show that for every unitary 
$d\times d$ matrix $u = (u_{ij})$, there is a $*$-automorphism 
$\theta_u$ of $\Cal T_d/\Cal K$ such that 
$$
\theta_u(Z_i) = \sum_{j=1}^d \bar u_{ji}Z_j.  
$$
For that, consider the unitary operator $U$ acting 
on $E$ by 
$$
Ue_i = \sum_{j=1}^d \bar u_{ji}e_j.  
$$
Then $\Gamma(U)$ is a unitary operator on $H^2$ for which 
$$
\Gamma(U)S_i\Gamma(U)^* = \sum_{j=1}^d\bar u_{ji}S_j, 
$$
hence $\theta_u$ is obtained by promoting the spatial 
automorphism $T\mapsto \Gamma(U)T\Gamma(U)^*$ of $\Cal T_d$ 
to the quotient $\Cal T_d/\Cal K$.  

The identification of $\Cal T_d/\Cal K$ with $C(\partial B_d)$ 
asserted by $\pi(S_i)=z_i$, $i=1,\dots,d$ is now obvious\qed

\enddemo


\subheading{6.  $d$-contractions and $\Cal A$-morphisms}

The purpose of this section is to make some observations 
about the role of $\Cal A$-morphisms in function theory and
operator theory.  

\proclaim{Definition 6.1}
Let $\Cal A$ be a subalgebra of a unital \cstar\ $\Cal B$ 
which contains the unit of $\Cal B$.  An $\Cal A$-morphism 
is a completely positive linear map $\phi: \Cal B\to\Cal B(H)$ 
of $\Cal B$ into the operators on a Hilbert space $H$ such 
that $\phi(\bold 1)=\bold 1$ and 
$$
\phi(AX) = \phi(A)\phi(X)\qquad A\in \Cal A, X\in \Cal B. 
$$
\endproclaim

$\Cal A$-morphisms arose naturally in our work on the 
dilation theory of completely positive maps and semigroups
\cite{7,8,9}.  Jim Agler has pointed out that they are 
related to his notion of hereditary polynomials and 
hereditary isomorphisms 
(for example, see \cite{1, Theorem 1.5}).  Indeed, if 
$\Cal B$ denotes the $C^*$-algebra generated by a single 
operator $T$ and the identity, then one can 
show that a completely positive map of $\Cal B$ 
which is a hereditary isomorphism
on the space of hereditary polynomials in $T$ 
is an $\Cal A$-morphism relative to the algebra $\Cal A$
of all polynomials in the {\it adjoint} $T^*$.  

In general the restriction of an $\Cal A$-morphism to 
$\Cal A$ is a completely contractive representation of 
the subalgebra $\Cal A$ on $H$.  Theorem 4.5 implies that
every $d$-contraction $\bar T$ acting on a Hilbert space $H$ 
gives rise to a contraction $L: \Cal F_+(\Bbb C^d)\otimes K\to H$
which intertwines the action of the $d$-shift and $\bar T$.  
$L$ is often a co-isometry, and that implies 
the following assertion about $\Cal A$-morphisms.   

\proclaim{Theorem 6.2}
Let $\Cal A$ be the subalgebra of the Toeplitz \cstar\ 
$\Cal T_d$ consisting of all polynomials 
in the $d$-shift $(S_1,\dots,S_d)$.  Then for every 
$d$-contraction $(T_1,\dots,T_d)$ acting on a Hilbert space
$H$ there is a unique 
$\Cal A$-morphism 
$$
\phi: \Cal T_d\to \Cal B(H)
$$
such that $\phi(S_k)=T_k$, $k=1,\dots,d$.  

Conversely, every $\Cal A$-morphism $\phi:\Cal T_d\to\Cal B(H)$
gives rise to a $d$-contraction $(T_1,\dots,T_d)$ on $H$ by way of 
$T_k=\phi(S_k)$, $k=1,\dots,d$.  
\endproclaim

\demo{proof}
The uniqueness assertion is immediate from 
(5.8), since an $\Cal A$-morphism is 
uniquely determined on the closed
linear span of the set of products $\{AB^*: A,B\in \Cal A\}$.  

For existence, we first show that every 
{\it null} $d$-contraction $\bar T=(T_1,\dots,T_d)$ defines 
an $\Cal A$-morphism as asserted in Theorem 6.2.  For that,
let 
$$
\Delta = (\bold 1-T_1T_1^*-\dots-T_dT_d^*)^{1/2},  
$$
let $K=\overline{\Delta H}$ be the closed range of 
$\Delta$ and let $\Cal F_+(E)$ be the symmetric Fock 
space over $E=\Bbb C^d$.  Choose an orthonormal basis 
$e_1,\dots,e_d$ for $E$.  Theorem 4.5 asserts that there 
is a unique bounded operator $L:\Cal F_+(E)\otimes K\to H$
satisfying $L(1\otimes \xi)=\Delta\xi$ for $\xi\in K$ and 
$$
L(e_{i_1}e_{i_2}\dots e_{i_n}\otimes \xi)=
T_{i_1}T_{i_2}\dots T_{i_n}\Delta \xi, \tag{6.3}
$$
for $n=1,2,\dots$, $i_1,i_2,\dots,i_n\in\{1,\dots,d\}$, 
$\xi\in K$; moreover, since $(T_1,\dots,T_d)$ is a 
null $d$-contraction $L$ is a coisometry.  

We may consider that the $d$-shift $(S_1,\dots,S_d)$
is defined on $\Cal F_+(E)$ by 
$$
S_k\xi = e_k\xi,\qquad k=1,\dots,d.  
$$
(6.3) implies that 
$$
L(f(S_1,\dots,S_d)\otimes \bold 1_K)=f(T_1,\dots,T_d)L 
\tag{6.4}
$$
for every polynomial $f$ in $d$ variables.  Let
$\phi: \Cal T_d\to\Cal B(H)$ be the completely positive map 
$$
\phi(X) = L(X\otimes\bold 1_K) L^*,\qquad X\in \Cal T_d.  
$$
Since $L^*$ is an isometry we have $\phi(\bold 1)=\bold 1_H$.  
(6.4) implies that for every $X\in\Cal T_d$ we have 
$$
\phi(f(S_1,\dots,S_d)X)=f(T_1,\dots,T_d)\phi(X), 
$$
hence $\phi$ is an $\Cal A$-morphism having the required properties.  

The general case is deduced from this by a simple device.  
Let $\bar T=(T_1,\dots,T_d)$ 
be any $d$-contraction, choose a number $r$ 
so that $0<r<1$, and set 
$$
\bar T_r=(rT_1,\dots,rT_d).  
$$
The row norm of the $d$-tuple $\bar T_r$ is at 
most $r$, hence $\bar T_r$ is a null $d$-contraction
(see Remark 4.4).  
By what was just proved there is an $\Cal A$-morphism 
$\phi_r: \Cal T_d\to \Cal B(H)$ satisfying 
$$
\phi_r(S_k)=rT_k,\qquad k=1,\dots,d.  
$$
We have 
$$
\phi_r(f(S_1,\dots,S_d)g(S_1,\dots,S_d)^*)=
f(rT_1,\dots,rT_d)g(rT_1,\dots,rT_d)^*
$$
for polynomials $f$, $g$.  
Since operators of the form 
$f(S_1,\dots,S_d)g(S_1,\dots,S_d)^*$ span $\Cal T_d$ and 
since the family of maps $\phi_r$, $0<r<1$ is uniformly 
bounded, it 
follows that $\phi_r$ converges
point-norm to an $\Cal A$-morphism $\phi$ as $r\uparrow 1$,
and $\phi(S_k)=T_k$ for all $k$.  

It remains only to show that for every $\Cal A$-morphism 
$\phi: \Cal T_d\to \Cal B(H)$, the operators 
$T_k=\phi(S_k)$ define a $d$-contraction.  To see that, 
write 
$$
T_kT_k^*=\phi(S_k)\phi(S_k)^* = \phi(S_kS_k^*).  
$$
Then 
$$
\sum_{k=1}^dT_kT_k^* = \phi(\sum_{k=1}^d S_kS_k^*) \leq 
\phi(\bold 1)=\bold 1.  
$$
So by remark 3.2, $(T_1,\dots,T_d)$ is a $d$-contraction\qed
\enddemo

\remark{Remarks}
We have already pointed out that in 
general, an $\Cal A$-morphism must be a completely 
contractive representation of $\Cal A$.  Conversely, if 
$\Cal A$ is the polynomial algebra in $\Cal T_d$ and 
$\phi: \Cal A\to\Cal B(H)$ is a representation which 
is $d$-contractive in the sense that its natural 
promotion to $d\times d$ matrices over $\Cal A$ is 
a contraction, then after noting that the operator 
matrix 
$$
A = 
\pmatrix
S_1&S_2&\hdots&S_d\\
0&0&\hdots&0\\
\vdots&\vdots &&\vdots\\
0&0&\hdots&0
\endpmatrix
\in M_d(\Cal T_d)
$$
satisfies $\|A\|^2=\|AA^*\| = \|S_1S_1^*+\dots+S_dS_d^*\|=1$, 
we find that the image of $A$ under the promotion 
of $\phi$ is a contraction, and hence $T_k=\phi(S_k)$, 
$k=1,\dots,d$ defines a $d$-contraction.  Thus, we may conclude

\proclaim{Corollary 1}Let $d=1,2,\dots$.  
Every $d$-contractive representation $\phi$ of the 
polynomial algebra $\Cal A\subseteq \Cal T_d$ 
is completely contractive, and 
can be extended uniquely
to an $\Cal A$-morphism 
$$
\tilde\phi:\Cal T_d\to\Cal B(H).  
$$
\endproclaim

We have already seen that the unitary group $\Cal U_d$ of 
$\Bbb C^d$ acts naturally on $\Cal T_d$ as a group of 
$*$-automorhisms by way of 
$$
\theta_U(X) = \Gamma(U)X\Gamma(U)^*,\qquad X\in\Cal T_d, U\in\Cal U_d.  
$$
As a straightforward application of Theorem 6.2 we show that 
the definition of $\theta$ can be extended to all contractions 
in $\Cal B(\Bbb C^d)$ so as to obtain a semigroup of 
$\Cal A$-morphisms acting on $\Cal T_d$.  

\proclaim{Corollary 2}Let $\Cal A\subseteq \Cal T_d$ 
be the algebra of all polynomials in $S_1,\dots, S_d$.  
For every contraction $A$ acting on $\Bbb C^d$ there is a unique 
$\Cal A$-morphism $\theta_A: \Cal T_d\to\Cal B(H^2)$ satisfying 
$$
\theta_A(M_f) = M_{f\circ A^*} \tag{6.5}
$$
for every linear functional $f$ on $\Bbb C^d$, $A^*$ denoting 
the adjoint of $A\in\Cal B(\Bbb C^d)$.    
\endproclaim
\demo{proof}
Considering the polar decomposition of $A$, we may find 
a pair of orthonormal bases 
$u_1,\dots,u_d$ and $u_1^\prime,\dots,u_d^\prime$ for 
$\Bbb C^d$ and numbers $\lambda_k$ in the unit interval 
such that 
$$
Au_k = \lambda_ku_k^\prime, \qquad k=1,\dots,d.  
$$
Let $z_1,\dots,z_d$ and $z_1^\prime,\dots,z_d^\prime$ be the 
corresponding systems of orthogonal coordinate functions
$$
\align
z_k(x)&=\<x,u_k\>,\\
z_k^\prime(x)&=\<x,u_k^\prime\>.  
\endalign
$$
The linear functionals $z_k,z_k^\prime$ are related 
by 
$$
z_k\circ A^* = \lambda_k z_k^\prime,\qquad k=1,\dots,d.  \tag{6.6}
$$
Thus if we realize the $d$-shift $(S_1,\dots,S_d)$ 
as $S_k=M_{z_k}$ and if we set
$T_k=\lambda_k M_{z_k^\prime}$, then $(T_1,\dots,T_d)$ is a 
$d$-contraction and Theorem 6.2 implies that there is a unique 
$\Cal A$-morphism $\theta_A:\Cal T_d\to\Cal B(H^2)$ such 
that $\theta_A(S_k)=T_k$ for every $k$.  After noting that 
$\theta_A$ satisfies (6.5) because of (6.6) above, the proof 
is complete\qed
\enddemo

From (6.5) together with the uniqueness assertion 
of Theorem 6.2 it follows 
that for two contractions $A,B\in\Cal B(\Bbb C^d)$ 
we have $\theta_{AB}=\theta_A\circ\theta_B$.  
It is routine to verify that 
$\theta_A(\Cal T_d)\subseteq \Cal T_d$, that for every fixed 
$X\in\Cal T_d$
the function $A\mapsto \theta_A(X)$ moves continuously in 
the norm of $\Cal T_d$, and that $\theta_A$ agrees with 
the previous definition when $A$ is unitary.  
\endremark

\remark{Uniqueness of Representing measures}
Representing measures for points in the interior of 
the unit ball in $\Bbb C^d$ are notoriously non-unique
in dimension $d\geq 2$.  Indeed, for every 
$\bar t=(t_1,\dots,t_d)\in B_d$ there is an 
uncountable family of probability measures 
$\mu_\alpha$ supported in the boundary $\partial B_d$
such that $\mu_\alpha\perp\mu_\beta$ 
for $\alpha\neq\beta$ and 
$$
\int_{\partial B_d}f(\zeta)\,d\mu_\alpha(\zeta) = 
f(\bar t), \qquad f\in\Cal P  
$$
see \cite{37, p. 186}.  The following result asserts that one 
can recover uniqueness by replacing
measures on $\partial B_d$ with states on the 
Toeplitz \cstar\ which define $\Cal A$-morphisms.
\endremark

\proclaim{Corollary 3}
Let $\bar t=(t_1,\dots,t_d)\in\Bbb C^d$ satisfy 
$|t_1|^2+\dots+|t_d|^2<1$ and let $\bar S=(S_1,\dots,S_d)$
be the $d$-shift.  Then there is a unique state 
$\phi$ of $\Cal T_d$ satisfying 
$$
\phi(f(\bar S)g(\bar S)^*)=f(\bar t)\bar g(\bar t), 
\qquad f,g\in \Cal P.  \tag{6.7}
$$
$\phi$ is the (pure) vector state 
$$
\phi(A) = (1-\|\bar t\|^2)\<Au_{\bar t},u_{\bar t}\>, 
\qquad A\in\Cal T_d 
$$
where $u_{\bar t}(x) = (1-\<x,\bar t\>)^{-1}$ is the 
$H^2$ function defined in (1.11).  
\endproclaim
\demo{proof}
We may consider that $\bar t=(t_1,\dots,t_d)$ is a 
$d$-contraction acting on the one-dimensional Hilbert 
space $\Bbb C$.  Theorem 6.2 implies that there is a 
unique state $\phi: \Cal T_d\to\Bbb C$ satisfying (6.7), 
and it remains only to identify $\phi$.  From (2.4) 
we have 
$$
\<M_fM_g^*u_{\bar t},u_{\bar t}\>=
\<M_g^*u_{\bar t},M_f^*u_{\bar t}\>=
f(\bar t)\bar g(\bar t)\|u_{\bar t}\|^2=
(1-\|\bar t\|^2)^{-1}f(\bar t)\bar g(\bar t)
$$
as asserted.\qed
\enddemo


\subheading{7.  The $d$-shift as an operator space}

In this section we consider the operator space 
$\Cal S_d\subseteq \Cal B(H^2)$ 
generated by the $d$-shift $(S_1,\dots,S_d)$,
$$
\Cal S_d = \{a_1S_1+\dots+a_dS_d: a_1,\dots,a_d\in \Bbb C\}.  
$$
By a commutative operator space we mean a linear subspace 
$\Cal S\subseteq \Cal B(H)$ whose operators mutually
commute with one another.  
We introduce a sequence of numerical invariants for 
arbitrary operator spaces, and for dimension $d\geq 2$ 
we show that among all $d$-dimensional commutative 
operator spaces, $\Cal S_d$ 
is distinguished by the fact that its sequence of 
numerical invariants is maximal (Theorem 7.7).   

Given an arbitrary operator space $\Cal S\subseteq\Cal B(H)$, 
let $\bar T=(T_1,T_2,\dots)$ be an infinite sequence 
of operators in $\Cal S$ such that all but a finite 
number of terms are $0$.  We write $\seq(\Cal S)$ 
for the set of all such sequences.  Every such sequence 
has a ``row norm" and a ``column norm", depending on whether
one thinks of the sequence as defining an operator 
in $\Cal B(H^\infty,H)$ or in $\Cal B(H,H^\infty)$.  
These two norms are familiar and easily computed,  
$$
\align
\|\bar T\|_{row} &= \|\sum_k T_kT_k^*\|^{1/2}\\
\|\bar T\|_{col} &= \|\sum_k T_k^*T_k\|^{1/2}.  
\endalign
$$
Given two sequences $\bar T, \bar T^\prime\in\seq (\Cal S)$,
we can form a product sequence $(T_iT_j^\prime: i,j=1,2,\dots)$
which we may consider an element of $\seq (\Cal B(H))$, if we 
wish, by relabelling the double sequence as a single sequence.  
Though for the computations below it will be more convenient 
to allow the index set to vary in the obvious way.  
In particular, every $\bar T\in\text{seq}(\Cal S)$ 
can be raised to the $n$th power to obtain 
$\bar T^n\in\seq (\Cal B(H))$, $n=1,2,\dots$.  
For each $n=1,2,\dots$ we define $E_n(\Cal S)\in[0,+\infty]$ 
as follows,
$$
E_n(\Cal S) = \sup\{\|\bar T^n\|_{col}^2: 
\bar T\in\seq (\Cal S), \|\bar T\|_{row}\leq 1\}.  
$$
In the most explicit terms, we have
$$
E_n(\Cal S) = 
\sup\{\|\sum_{i_1,\dots,i_n = 1}^\infty T_{i_1}^*\dots 
T_{i_n}^*T_{i_n}\dots T_{i_1}\|: T_i\in \Cal S, 
\|\sum_{i=1}^\infty T_iT_i^*\|\leq 1 \},  
$$
the sup being taken over finitely nonzero sequences
$T_i\in\Cal S$.  

\proclaim{Definition 7.1}
$E_1(\Cal S), E_2(\Cal S),\dots$ is called the energy sequence 
of the operator space $\Cal S$.  
\endproclaim

If $\Cal S$ is the 
one-dimensional space spanned by a single operator
$T$ of norm $1$, then the energy sequence degenerates 
to $E_n(\Cal S) = \|T^n\|^2$, $n=1,2,\dots$.  In general, 
$E_n(\Cal S)^{1/2}$ is the norm of the homogeneous polynomial 
$\bar T\mapsto \bar T^n$, considered as a map of 
row sequences in $\Cal S$ to column sequences in $\Cal B(H)$.

\remark{Remarks}
We have defined the energy sequence 
in elementary terms.  It is useful,
however, to relate it to completely positive maps.  
Fixing an operator space $\Cal S$, notice that every 
sequence $\bar T\in\seq (\Cal S)$ gives rise to a normal completely
positive map $P_{\bar T}$ on $\Cal B(H)$ as the sum of the 
{\it finite} series 
$$
P_{\bar T}(A)=T_1AT_1^*+T_2AT_2^*+\dots . \tag{7.2}  
$$
Let $\cp (\Cal S)$ denote the set of all completely 
positive maps of the form (7.2).  The norm of $P = P_{\bar T}$ 
is given by 
$$
\|P\| = \|P(\bold 1)\| = \|\bar T\|_{row}.  
$$
Now any map $P\in \cp (\Cal S)$ of the form (7.2) has 
an adjoint $P_*$ which is defined as the completely positive 
map satisfying 
$$
\tr(P(A)B) = \tr(AP_*(B))
$$
for all finite rank operators $A,B$.  One finds
that if $P \in\cp(\Cal S)$ is given by the finitely nonzero sequence 
$\bar T$ then $P_*\in\cp (\Cal S^*)$ is 
given by the sequence of adjoints 
$$
P_*(A) = T_1^*AT_1+T_2^*AT_2+\dots.  \tag{7.3}
$$
Of course $P$, being a normal linear map of $\Cal B(H)$, 
is the adjoint of a bounded linear map $P_*$ acting on 
the predual of $\Cal B(H)$, and the map of (7.3) is simply
this preadjoint extended from the trace class operators 
to all of $\Cal B(H)$ (note that we use the fact that the 
sequence $\bar T$ is finitely nonzero here, since in general 
a bounded linear map of the trace class operators can be
unbounded relative to the operator norm, and thus not extendable
up to $\Cal B(H)$).  
\endremark
In any case, we find that if $P\in\cp (\Cal S)$ has 
the form $P = P_{\bar T}$ for $\bar T\in\seq (\Cal S)$ then 
$$
\|P_*\| = \|P_*(\bold 1)\| = \|\bar T\|_{col}.  
$$
Thus the definition of $E_n(\Cal S)$ can be restated as 
follows,
$$
E_n(\Cal S) = \sup\{\|P_*^n\| : P\in\cp (\Cal S), \|P\|\leq 1\}.  
\tag{7.4}
$$
The following result implies that for a finite dimensional 
operator space $\Cal S$ the terms of the energy sequence 
are all finite, and if $\Cal S$ is {\it commutative} then 
they grow no faster than $E_n(\Cal S)=O(n^{d-1})$, where 
$d$ is the dimension of $\Cal S$.  

\proclaim{Proposition 7.5}Let $\Cal S$ be an operator space 
of finite dimension $d$.  Then 
$$
E_n(\Cal S)\leq d^n, 
$$  
and if $\Cal S$ is also commutative then 
$$
E_n(\Cal S)\leq \frac{d(d+1)\dots (d+n-1)}{n!} = 
\frac{(n+d-1)!}{n!(d-1)!}.  
$$

\endproclaim

\demo{proof}
Let $P\in \cp(\Cal S)$ satisfy $\|P\|\leq 1$, and
let $d$ be the dimension of $\Cal S$.  It is clear that the 
metric operator space \cite{8} of $P$ is a subspace of 
$\Cal S$, and in particular there is a linearly independent 
set of $r\leq d$ elements $T_1,\dots,T_r$ in $\Cal S$ such that 
$$
P(A) = T_1AT_1^*+\dots T_rAT_r^*, \qquad A\in \Cal B(H).  
$$
Since $\|P\|=\|P(\bold 1)\| = \|T_1T_1^*+\dots+T_rT_r^*\|\leq 1$ 
it follows that $\|T_k\|\leq 1$ for every $k=1,\dots,r$, hence
$$
\|P_*\|=\|T_1^*T_1+\dots+T_r^*T_r\|\leq r\leq d.  
$$
Thus $\|P_*^n\|\leq d^n$ for every $n=1,2,\dots$.  From 
(7.4) we conclude that $E_n(\Cal S)\leq d^n$.  

In fact, the preceding argument shows that if 
$Q$ is a normal completely positive map of $\Cal B(H)$ 
whose metric operator space is $r$-dimensional and which 
satisfies $\|Q\|\leq 1$, then we have $\|Q_*\|\leq r$.  

We apply this to $Q=P^n$ as follows.  
By \cite{8}, the metric operator 
space $\Cal E_n$ of $P^n$ is a subspace of 
$$
\text{span}\{L_1L_2\dots L_n: L_i\in \Cal S\}.  
$$
Assuming $\Cal S$ to be commutative, the latter is naturally 
isomorphic to a quotient of the $n$-fold {\it symmetric} 
tensor product of vector spaces $\Cal S^n$.  Since 
$$
\dim \Cal S^n = \frac{(n+d-1)!}{n!(d-1)!}
$$
(see formula A.5 of the Appendix), and since $\|P^n\|\leq 1$, 
we find that 
$$
\|P_*^n\| = \|(P^n)_*\| \leq \dim \Cal E_n\leq 
\frac{(n+d-1)!}{n!(d-1)!}.
$$
The required estimate follows from the observation (7.4).\qed
\enddemo

\remark{Remark}
The asserted growth rate of the binomial coefficients of 
Proposition 7.5 is well known, and the precise asymptotic
relation is reiterated in formula A.6.  
\endremark

Throughout the remainder of this section we will be concerned
with finite dimensional commutative operator spaces.  

\proclaim{Definition 7.6}
A commutative operator space $\Cal S$ of finite 
dimension $d$ is said to be maximal if for every 
$n=1,2,\dots$ we have 
$$
E_n(\Cal S) = \frac{(n+d-1)!}{n!(d-1)!}.  
$$
\endproclaim

\remark{Remarks}
It is obvious that if $\Cal S$ is a space 
of mutually commuting normal operators, then 
$E_n(\Cal S)= 1$ for every $n$.  Similarly, it is 
not hard to show that if $\Cal S$ is a 
finite dimensional space of 
commuting quasinilpotent operators, then 
$$
\lim_{n\to\infty}E_n(\Cal S) = 0.  
$$
Thus the maximal spaces are 
rather far removed from both of these types.  

It is also true (though less obvious) that if $\Cal S$ 
is a commutative operator space of dimension $d$ for which 
$$
E_n(\Cal S) =\frac{(n+d-1)!}{n!(d-1)!}
$$
for some particular value of $n\geq 2$, then 
$$
E_k(\Cal S) =\frac{(k+d-1)!}{k!(d-1)!}
$$
for every $k=1,2,\dots,n$.  Thus for operator 
spaces which are not maximal, once the sequence of 
numbers $E_n(\Cal S)$ departs from the sequence 
of maximum possible values, it never returns.  
We omit the proof of these assertions since 
they are not required in the sequel.   
\endremark

\proclaim{Theorem 7.7}
For every $d=1,2,\dots$ the operator space $\Cal S_d$ of
the $d$-shift is maximal. 

Conversely, if $d\geq 2$ and if $\Cal S$ is a 
$d$-dimensional commutative operator space which 
is maximal, then there is a 
representation 
$\pi$ of the unital \cstar\ $C^*(\Cal S)$ 
generated by $\Cal S$ on $H^2$ such that

\roster
\item
$\pi(\Cal S) = \Cal S_d$, and 
\item
the restriction of $\pi$ to the unital subalgebra 
of $C^*(\Cal S)$ generated by $\Cal S$ is completely isometric.  
\endroster
In particular, the Toeplitz \cstar\ $\Cal T_d$ is 
isomorphic to a quotient of $C^*(\Cal S)$.  
\endproclaim

Before giving the proof of 7.7, we deduce from it the following 
characterization of $\Cal S_d$ as a space of 
essentially normal operators (by that we mean a 
commuting family of operators 
in $\Cal B(H)$ whose image 
in the Calkin algebra consists of normal elements).  
We remark that both the Corollary
and the essential part of Theorem 7.7 are false 
in dimension one.  

\proclaim{Corollary} Assume $d\geq 2$.  
Up to unitary equivalence, the space $\Cal S_d$ spanned 
by the $d$-shift is the only $d$-dimensional irreducible 
commutative operator space consisting of 
essentially normal operators, which is maximal 
in the sense of Definition 7.6.  
\endproclaim

\demo{proof of Corollary}
Suppose that $\Cal S$ acts on a Hilbert space $H$, 
and let $\Cal K$ denote the algebra of all compact 
operators on $H$.  
Let $\pi: C^*(\Cal S)\to\Cal B(H^2)$ be the representation
of Theorem 7.7.  The operators in $\Cal S$ cannot be 
normal because $\Cal S_d=\pi(\Cal S)$ contains no normal operators.  
Since 
$[\Cal S^*,\Cal S]\subseteq \Cal K\cap C^*(\Cal S)$
and since $C^*(\Cal S)$ is irreducible, it follows 
that $C^*(\Cal S)$ contains $\Cal K$.  $\pi(\Cal K)$ 
cannot be $\{0\}$ because that would imply that 
$\pi(\Cal S)=\Cal S_d$ consists of normal operators.  

Thus $\pi$ is an irreducible representation of 
$C^*(\Cal S)$ which is nonzero on $\Cal K$, hence
$\pi$ must be unitarily equivalent to the identity 
representation of $C^*(\Cal S)$.  
In particular, $\Cal S$ is 
unitarily equivalent to $\Cal S_d$\qed
\enddemo

\demo{proof of Theorem 7.7}  The proof of Theorem 
7.7 will occupy the remainder of this section.  
Let $(S_1,\dots,S_d)$ be the $d$-shift, let 
$\Cal S_d=\text{span}\{S_1,\dots,S_d\}$ be its associated 
operator space, and define $P\in\cp(\Cal S_d)$ by
$$
P(A) = S_1AS_1^*+\dots+S_dAS_d^*.  
$$
By the remark following Definition 2.10 
we have $P(\bold 1) = \bold 1-E_0$, hence
$\|P\|=1$.  Thus to show that $\Cal S_d$ is maximal
it suffices to show that for each 
$n\geq 1$, the operator $P_*^n(\bold 1)$ satisfies
$$
\|P_*^n\|=\|P_*^n(\bold 1)\| = \frac{(n+d-1)!}{n!(d-1)!}.  \tag{7.8}
$$
While (7.8) can be deduced directly from Corollary 4.9, 
we actually require somewhat more information 
about the operators $P_*^n(\bold 1)$ and their 
eigenvalue distributions.  

\proclaim{Lemma 7.9}
Let $N=E_1+2E_2+3E_3+\dots$ 
be the number operator acting on $H^2$ and for every 
$n=1,2,\dots$ let $g_n:[0,\infty)\to\Bbb R$ be the bounded
continuous function 
$$
g_n(x) = \prod_{k=1}^{n}\frac{x+k+d-1}{x+k}.  
$$
Then 
$$
P_*^n(\bold 1) = g_n(N) = \sum_{k=0}^\infty g_n(k)E_k.  
$$
The eigenvalue sequence $\{g_n(0)\geq g_n(1)\geq\dots\}$ of 
$P_*^n(\bold 1)$ is decreasing and we have 
$$
\|P_*^n(\bold 1)\| = g_n(0) = \frac{(n+d-1)!}{n!(d-1)!}.  
$$
If $d\geq 2$ then the eigenvalue sequence is strictly decreasing,
$g_n(0)>g_n(1)>\dots$.  
\endproclaim
\demo{proof of Lemma 7.9}
The assertions follow from a direct computation, which 
can be organized as follows.  By Proposition 5.3 we have 
$$
P_*(\bold 1) = g_1(N), \tag{7.10}
$$
where $N$ is the number operator and $g_1$ is the function 
of a real variable defined by 
$$
g_1(x) = \frac{x+d}{x+1}, \qquad x\geq 0.  
$$
More generally, if $g$ is any bounded continuous function 
defined on $[0,\infty)$, then we have 
$$
P_*(g(N)) = \tilde g(N) \tag{7.11}
$$
where 
$$
\tilde g(x) = g(x+1)\frac{x+d}{x+1}, \qquad x\geq 0.  
$$
Indeed, (7.11) follows from the fact that if 
$E_k$ denotes the $k$th spectral projection of $N$,
$$
N = \sum_{k=1}^\infty k E_k, 
$$
then $E_k$ is the projection on the subspace of homogeneous
polynomials of degree $k$ in $H^2$, and thus for each $i=1,\dots,d$
we have the commutation formulas 
$S_i^*E_0 = 0$,  and $S_i^*E_k=E_{k-1}S_i^*$ for 
$k\geq 1$.  It follows that $P_*(E_0)=0$ and 
$P_*(E_k)=E_{k-1}P_*(\bold 1)$ for $k=1,2,\dots$, thus
$$
P_*(g(N))=\sum_{k=1}^\infty g(k)P_*(E_k)=\tilde g(N).  
$$

After iterating (7.11) we find that 
$P_*^n(\bold 1) = g_n(N)$ where 
$$
g_n(x) = \prod_{k=1}^{n}\frac{x+d-1}{x+k}.  
$$
Since each $g_n$ is a monotone decreasing function we conclude 
that 
$$
\|P_*^n(\bold 1)\| = g_n(0) = \frac{d(d+1)\dots (d+n-1)}{n!},
$$
and (7.8) follows.  It is clear from the recurrence 
formula for $g_{n+1}$ in terms of $g_n$ that when $d\geq 2$,  
$g_n(x)$ is strictly decreasing in $x$\qed
\enddemo

\proclaim{Corollary}
Let $\omega$ be a state of the Toeplitz algebra 
$\Cal T_d$, $d\geq 2$, such that for some $n\geq 1$ 
we have 
$$
\omega(P_*^n(\bold 1)) = \|P_*^n(\bold 1)\| = 
\frac{(n+d-1)!}{n!(d-1)!}.  
$$
Then $\omega$ is the ground state $\omega(X)=\<Xv,v\>$, 
$v$ denoting the constant polynomial $v=1$.  
\endproclaim
\demo{proof}
Fix $n$.  By Lemma 7.9 we have 
$$
P_*^n(\bold 1) = \lambda_0 E_0+\lambda_1 E_1+\dots
$$
where $\lambda_0>\lambda_1>\dots>0$ and 
$\lambda_0=\|P_*^n(\bold 1)\|$.  Thus $P_*^n(\bold 1)$ has 
the form 
$$
P_*^n(\bold 1) = \lambda_0(E_0+K),
$$
where $K$ is a positive operator satisfying 
$K=(\bold 1-E_0)K(\bold 1-E_0)$ and $\|K\|=\lambda_1/\lambda_0<1$.  
Since $\omega(P_*^n(\bold 1))=\lambda_0$ we have 
$$
\omega(E_0)+\omega(K) = 1.  
$$
If $\omega(E_0)<1$ then we would have
$$
\omega(K)\leq \|K\|\omega(\bold 1-E_0)=
\|K\|( 1-\omega(E_0))<1-\omega(E_0), 
$$
contradicting the preceding equation.  Hence $\omega(E_0)=1$
and $\omega$ must be the ground state\qed
\enddemo

In particular, Lemma 7.9 implies 
that $\Cal S_d$ is maximal among all $d$-dimensional 
commutative operator spaces.  

In order to prove the converse assertion of Theorem 7.7, 
we recall one or two facts 
from the theory of boundary representations (see
\cite{3}, 2.1.2 and 2.2.2).  By a unital 
operator space we mean a pair $\Cal S\subseteq \Cal B$ 
consisting of a linear subspace $\Cal S$ of a unital 
\cstar\ $\Cal B$, which contains the unit of $\Cal B$ and 
generates $\Cal B$ as a \cstar, $\Cal B = C^*(\Cal S)$.  
An irreducible representation $\pi: \Cal B\to \Cal B(H)$
is said to be a boundary representation for $\Cal S$ if
$\pi\restriction_{\Cal S}$ has a {\it unique} completely 
positive linear extension to $\Cal B$, namely $\pi$ 
itself.  Boundary representations are the noncommutative 
counterpart of points in the Choquet boundary of a 
function space $S\subseteq C(X)$.  
Their key property 
is their functoriality; if $\Cal S_1\subseteq \Cal B_1$
and $\Cal S_2\subseteq\Cal B_2$ are unital operator spaces and 
$\phi: \Cal S_1\to \Cal S_2$ is a completely isometric linear 
map satisfying $\phi(\bold 1)=\bold 1$ and 
$\phi(\Cal S_1)=\Cal S_2$, then for every boundary representation
$\pi_2: \Cal B_2\to\Cal B(H)$ for $\Cal S_2$ there is a unique
boundary representation $\pi_1: \Cal B_1\to\Cal B(H)$ for 
$\Cal S_1$ which satisfies 
$$
\pi_2(\phi(T)) = \pi_1(T), \qquad T\in \Cal S_1.  \tag{7.12}
$$

\proclaim{Lemma 7.13}
For $d\geq 2$, the identity representation of the 
Toeplitz algebra $\Cal T_d$ is a boundary representation 
for the $d+1$-dimensional subspace 
$\text{span}\{\bold 1,S_1,\dots,S_d\}$.  
\endproclaim
\demo{proof}
By \cite{4, Theorem 2.1.1} it is enough to show that 
the Calkin map is not isometric when promoted to the 
space $M_d\otimes \Cal S$ of $d\times d$ matrices over $\Cal S$. 
Consider the operator $A\in M_d\otimes \Cal S$  defined by
$$
A=
\pmatrix
S_1&0&\hdots&0\\
S_2&0&\hdots&0\\
\vdots&\vdots &&\vdots\\
S_d&0&\hdots&0
\endpmatrix
$$
Then $\|A^*A\|=\|S_1^*S_1+\dots +S_d^*S_d\|=d$ 
by Proposition 5.3.  Hence
$\|A\|=\sqrt d$.  On the other hand, 
by Theorem 5.7 
the Calkin map carries $S_k$ to the $k$th coordinate 
function $z_k(x)=\<x,e_k\>$, $x\in\partial B_d$.   
Hence the image of $A$ under the promoted 
Calkin map is the matrix of functions
on $\partial B_d$ defined by
$$
F(x)=
\pmatrix
z_1(x)&0&\hdots&0\\
z_2(x)&0&\hdots&0\\
\vdots&\vdots &&\vdots\\
z_d(x)&0&\hdots&0
\endpmatrix.  
$$
Clearly $\sup\{\|F(x)\|: x\in\partial B_d\}=1< \sqrt d=\|A\|$, 
as required.  
\qed
\enddemo

\proclaim{Lemma 7.14}
Let $\Cal S\subseteq \Cal B(H)$ be a commutative 
operator space of finite dimension $d\geq 2$, and suppose
$\Cal S$ is maximal.  Then there is a state $\rho$ of 
the unital \cstar\ $C^*(\Cal S)$ generated by $\Cal S$ 
and a $d$-contraction $\bar T=(T_1,\dots,T_d)$, $T_i\in\Cal S$, 
such that 
$$
\rho(g(\bar T)^*f(\bar T))=\<f,g\>_{H^2}, 
$$
for all polynomials $f,g\in\Cal P$.  
\endproclaim
\demo{proof of Lemma 7.14}
The set of all $d$-contractions $(T_1,\dots,T_d)$ 
whose component operators belong to $\Cal S$ can be regarded 
as a compact subset of the cartesian product of $d$ copies 
of the unit ball of $\Cal S$, and of course the state space 
of $C^*(\Cal S)$ is weak$^*$-compact.  Thus, after a routine
compactness argument (which we omit), the proof of Lemma 
7.14 reduces to establishing the following assertion: for 
every $n=1,2,\dots$ there is a pair $(\rho,\bar T)$ 
consisting of a state $\rho$ of 
$C^*(\Cal S)$ and a $d$-contraction $\bar T=(T_1,\dots,T_d)$
whose components belong to $\Cal S$ such that 
$$
\rho(g(\bar T)^*f(\bar T))=\<f,g\>_{H^2} \tag{7.15}
$$
for all polynomials $f,g\in \Cal P$ of degree $\leq n$.  

To prove the latter, since $E_n(\Cal S)=\frac{(n+d-1)!}{n!(d-1)!}$
we may find a completely positive map $P\in\cp(\Cal S)$ such 
that $\|P\|\leq 1$ and 
$$
\|P_*^n(\bold 1)\|=\frac{(n+d-1)!}{n!(d-1)!}  \tag{7.16}
$$
(note that the supremum of (7.4) is achieved here because 
the space $\{P\in\cp(\Cal S): \|P\|\leq 1\}$ is compact).  
Considering that the metric operator space of $P$ is a subspace 
of $\Cal S$ \cite{8} we can find a (linearly independent) 
set $T_1,\dots,T_r\in\Cal S$ such that 
$$
P(A) = T_1AT_1^*+\dots T_rAT_r^*, \qquad A\in\Cal B(H).  
$$
By appending $T_{r+1}=\dots=T_d=0$ to the sequence if necessary, 
we can assume that $r=d$.  Because 
$$
\|P\|=\|P(\bold 1)\| = \|T_1T_1^*+\dots+T_dT_d^*\|\leq 1, 
$$
$\bar T=(T_1,\dots,T_d)$ is a $d$-contraction for which 
(7.16) holds.  

Let $\rho$ be any state of $C^*(\Cal S)$ satisfying 
$$
\rho(P_*^n(\bold 1))=\|P_*^n(\bold 1)\| = \frac{(n+d-1)!}{n!(d-1)!},
$$
and consider the positive semidefinite inner product defined 
on $\Cal P$ by 
$$
\<f,g\> = \rho(g(\bar T)^*f(\bar T)).  \tag{7.17}
$$
One sees (after consideration of the GNS construction for the 
state $\rho$) that since $\bar T$ is a $d$-contraction, the 
Hilbert seminorm $\|f\|^2=\rho(f(\bar T)^*f(\bar T))$ satisfies
$$
\|z_1f_1+\dots+z_df_d\|^2\leq \|f_1\|^2+\dots+\|f_d\|^2
$$
for all polynomials $f_1,\dots,f_d\in\Cal P$. By Proposition 
4.2, $\|\cdot\|$ will be a contractive seminorm provided 
that 
$$
1\perp z_1\Cal P+\dots+z_d\Cal P
$$ 
in its associated inner product space; or equivalently, that
$$
\rho(T_kf(\bar T))=0, \qquad k=1,\dots,d, f\in \Cal P.  
$$
Since $\rho$ is a state, the latter will follow if we establish
$$
\rho(T_kT_k^*)=0, \qquad k=1,\dots,d.  \tag{7.18}
$$

To prove 7.18, let $\phi: \Cal T_d\to \Cal B(H)$ be 
an $\Cal A$-morhpism satisfying $\phi(S_k)=T_k$, $k=1,\dots,d$
(see Theorem 6.2), and let $\omega$ be the state of 
$\Cal T_d$ defined by $\omega=\rho\circ\phi$.  We claim 
that $\omega$ is the ground state of $\Cal T_d$.  Indeed, 
for every $n$-tuple of integers $1\leq i_1,\dots,i_n\leq d$ 
we have by the Schwarz inequality
$$
\phi(S_{i_1}^*\dots S_{i_n}^*S_{i_n}\dots S_{i_1}) \geq 
\phi(S_{i_n}\dots S_{i_1})^*\phi(S_{i_n}\dots S_{i_1}) =
T_{i_1}^*\dots T_{i_n}^*T_{i_n}\dots T_{i_1},
$$
and hence 
$$
\omega(S_{i_1}^*\dots S_{i_n}^*S_{i_n}\dots S_{i_1})\geq 
\rho(T_{i_1}^*\dots T_{i_n}^*T_{i_n}\dots T_{i_1}).  
$$
Summing over all such $n$-tuples we obtain 
$$
\sum_{i_1,\dots,i_n=1}^d
\omega(S_{i_1}^*\dots S_{i_n}^*S_{i_n}\dots S_{i_1}) \geq
\rho(P_*^n(\bold 1))=\frac{(n+d-1)!}{n!(d-1)!}.  
$$
The Corollary of Lemma 7.9 implies that $\omega$ 
must be the ground state of $\Cal T_d$.  
In particular, for each $k=1,\dots,d$ we have 
$$
\rho(T_kT_k^*) = \rho(\phi(S_k)\phi(S_k)^*)=\rho(\phi(S_kS_k^*))
=\omega(S_kS_k^*)=\|S_k^*1\|^2=0, 
$$
and (7.18) follows.  

It is clear that the Hilbert seminorm of (7.17) is normalized
so that $\|1\|^2=\rho(\bold 1)=1$, so by Theorem 4.3 we have 
$\|f\|\leq \|f\|_{H^2}$ for every $f\in \Cal P$.  Corollary
4.9 now implies that (7.15) is satisfied, and the proof
is complete.\qed
\enddemo

To complete the proof of Theorem 7.7, we find a 
$d$-contraction $\bar T=(T_1,\dots,T_d)$ and a state 
$\rho$ of $C^*(\Cal S)$ satisfying the conditions of Lemma 7.14.  
Let $\Cal A\subseteq \Cal T_d$ be the commutative 
algebra generated by the $d+1$ operators 
$\{\bold 1,S_1,\dots,S_d\}$.  From Theorem 6.2 there is 
an $\Cal A$-morphism $\phi: \Cal T_d\to\Cal B(H)$ satisfying 
$\phi(S_k)=T_k$ for $k=1,\dots,d$.  

The restriction of $\phi$ to 
$\Cal A$ is completely contractive, and we claim that it is 
completely isometric.  To prove that, it is enough exhibit 
a representation $\sigma$ of $C^*(\Cal S)$ on another Hilbert 
space $K$ and a subspace $K_0\subseteq K$ which is invariant 
under the set of operators $\sigma(T_1),\dots,\sigma(T_k)$ such 
that the $d$-tuple of restrictions 
$$
(\sigma(T_1)\restriction_{K_0},\dots,\sigma(T_d)\restriction_{K_0})
$$
is unitarily equivalent to the $d$-shift $(S_1,\dots,S_d)$.  
The latter follows immediately from the properties 
of $\rho$.  Indeed, 
the GNS construction applied to $\rho$ gives a representation 
$\sigma$ of $C^*(\Cal S)$ on a Hilbert space $K$ and a unit 
cyclic vector $\xi\in K$ such that 
$$
\rho(X)=\<\sigma(X)\xi,\xi\>, \qquad X\in C^*(\Cal S).  
$$
Letting $K_0$ be the closure of $\{\sigma(f(\bar T))\xi:f\in \Cal P\}$
we see that $K_0$ is an invariant subspace which, 
since $\rho(g(\bar T)^*f(\bar T))=\<f,g\>_{H^2}$ for 
all $f,g\in \Cal P$, has the stated properties.  

In particular, the dimension of the subspace of $\Cal S$ 
spanned by $T_1,\dots,T_d$ must be $d=\dim(\Cal S_d)$, and 
thus $T_1,\dots,T_d$ is a basis for $\Cal S$.  
Using Lemma 7.13 and the functoriality property (7.12) of 
boundary representations it follows that there is a unique 
represenation $\pi: C^*(\Cal S)\to\Cal B(H^2)$, which is a 
boundary represenation for the unital algebra of all polynomials
in $T_1,\dots,T_d$, such that $\pi(T_k)=S_d$, $k=1,\dots,d$.  
The representation $\pi$ obviously has the properties asserted 
in Theorem 7.7\qed

\enddemo


\subheading{8.  Various applications}

In this section we give several applications of 
the preceding results to function theory and 
multivariable operator theory.  These are a version 
of von Neumann's inequality for arbitrary $d$-contractions, 
a model theory for $d$-contractions based on the 
$d$-shift, and a discussion of the absence of inner 
functions in the multiplier algebra of the $d$-shift.  

We point out that Popescu  has established various versions
of von Neumann's inequality for non-commutative $d$-tuples of 
operators \cite{30,32,34,35}.  
Here, on the other hand, we are concerned 
with $d$-contractions.  The version of 
von Neumann's inequality that is appropriate 
for $d$-contractions is the following.  

\proclaim{Theorem 8.1}
Let $\bar T=(T_1,\dots,T_d)$ be an arbitrary $d$-contraction 
acting on a Hilbert space $H$.  Then for every polynomial 
$f$ in $d$ complex variables we have 
$$
\|f(T_1,\dots,T_d)\| \leq \|f\|_\Cal M,   
$$
$\|f\|_\Cal M$ being the norm of $f$ in the multiplier 
algebra $\Cal M$ of $H^2$.  

More generally, let $(S_1,\dots,S_d)$ be the $d$-shift 
and let 
$\Cal A\subseteq\Cal T_d$ be the algebra 
of all polynomials in $S_1,\dots,S_d$.  Then the representation
$f\mapsto f(T_1,\dots,T_d)$ defines a completely contractive 
representation of $\Cal A$.  
\endproclaim
\demo{proof}
The assertions are immediate consequences of Theorem 6.2, once 
one observes that $\|f\|_\Cal M = \|f(S_1,\dots,S_d)\|$.\qed
\enddemo

Turning now to models, we first
recall some of the literature of dilation theory in 
$d$-dimensions.  There are a number of positive results concerning
non-commutative models for non-commuting $d$-tuples which 
satisfy the conditions of remark 3.2.  
The first results along these lines 
are due to Frazho \cite{21} for pairs of operators.  Frazho's 
results were generalized by Bunce \cite{15} to 
$d$-tuples.  Popescu has clarified 
that work by showing that such a $d$-tuple can 
often be obtained by compressing a certain natural $d$-tuple 
of isometries acting on the full Fock space 
$\Cal F(\Bbb C^d)$ over $\Bbb C^d$ (the left creation operators) 
to a co-invariant subspace of $\Cal F(\Bbb C^d)$, and he has 
worked out a functional calculus for that situation 
\cite{28,29,30,31}.  We also point out some recent work of 
Davidson and Pitts 
\cite{18,19}, relating to the operator algebra generated
by the left creation operators on the full Fock space.  

There is relatively little in the literature of operator theory, 
however, that relates to uniqueness of dilations in 
higher dimensions (however, see \cite{11}).  Indeed, normal 
dilations for $d$-contractions, when they exist, are almost never 
unique.  On the other hand, recent results 
in the theory of semigroups of completely positive maps
do include uniqueness.  Generalizing work 
of Parathasarathy, B. V. R. Bhat \cite{14} 
has shown that a unital semigroup of completely positive maps
of a von Neumann algebra $M$ can be dilated uniquely 
to an $E_0$-semigroup 
acting on a larger von Neumann algebra $N$ which 
contains $M$ as 
a hereditary subalgebra.  A similar (and simpler) result holds
for single unital completely positive maps: there is a unique dilation 
to a unital endomorphism acting on a larger von Neumann algebra as 
above.  In the case where $M=\Cal B(H)$, the latter dilation 
theorem is closely related to the Bunce-Frazho theory of 
$d$-tuples by way of the metric operator space associated 
with a normal completely postive map of $\Cal B(H)$ 
\cite{8,9}.  SeLegue \cite{42} has succeeded in unifying 
these results.  
 
In the following discussion, we reformulate Theorem 6.2 
as a concrete assertion about $d$-contractions 
which parallels some of the principal 
assertions of the Sz.-Nagy Foias model 
theory of $1$-contractions \cite{43}. 
Much of Theorem 8.5 follows directly from Theorem 6.2 and 
standard lore on the representation theory of \cstar s.  
For completeness, we have given a full sketch of the argument.   

We recall some elementary facts about the 
representation theory of \cstar s such as $\Cal T_d$.  
Let $\pi:\Cal T_d\to\Cal B(H)$ be a nondegenerate 
$*$-representation of $\Cal T_d$ on a separable Hilbert space $H$.  
Because of the exact sequence of Theorem 5.7, standard 
results about the representations of the \cstar\ of 
compact operators imply 
that $\pi$ decomposes 
into a direct sum $\pi_1\oplus\pi_2$ where $\pi_1$ is a 
multiple of $n=0,1,2,\dots,\infty$ 
copies of the identity representation
of $\Cal T_d$ and $\pi_2$ is a representation which annihilates
$\Cal K$.  $\pi_1$ and $\pi_2$ are disjoint as representations
of $\Cal T_d$.  This decomposition is unique in 
the sense that if $\pi_1^\prime$ 
is another multiple of $n^\prime$ copies of the 
identity representation of $\Cal T_d$ and $\pi_2^\prime$
annihilates $\Cal K$, and if $\pi_1^\prime\oplus\pi_2^\prime$ 
is unitarily equivalent
to $\pi_1\oplus\pi_2$, then $n^\prime=n$ and 
$\pi_2^\prime$ is unitarily equivalent to $\pi_2$ 
\cite{5}.  

We will make use of these observations in a form that 
relates more directly to operator theory. 
 
\proclaim{Definition 8.2}
Let $d=1,2,\dots$.  By a spherical operator (of dimension $d$) 
we mean a $d$-tuple $(Z_1,\dots,Z_d)$ of commuting normal operators 
acting on a common Hilbert space such that 
$$
Z_1^*Z_1+\dots +Z_d^*Z_d=\bold 1.  
$$ 
\endproclaim

Spherical operators are the higher dimensional counterparts
of unitary operators.  For every spherical operator $(Z_1,\dots,Z_d)$
acting on $H$ there is a unique unital $*$-representation 
$\pi: C(\partial B_d)\to\Cal B(H)$ which carries the $d$-tuple 
of canonical coordinate functions to $(Z_1,\dots,Z_d)$.  This
relation between $d$-dimensional spherical operators and 
nondegenerate representations of $C(\partial B_d)$ is bijective.  

If $\bar T=(T_1,\dots,T_d)$ 
is an arbitrary $d$-tuple of operators acting on a common 
Hilbert space $H$ and $n$ is a nonnegative integer or $+\infty$ 
we will write $n\cdot\bar T= (n\cdot T_1,\dots,n\cdot T_d)$ 
for the $d$-tuple of operators 
acting on the direct sum of $n$ copies of $H$ defined by
$$
n\cdot T_k = \underbrace{T_k\oplus T_k\oplus\dots}_{\text{$n$ times}}, 
$$
where for $n=0$ the left side is interpreted the {\it nil}
operator, that is, no operator at all.  The direct sum 
of two $d$-tuples of operators is defined in the obvious 
way as a $d$-tuple acting on the direct sum of Hilbert spaces.  
The preceding remarks are summarized as follows.  

\proclaim{Proposition 8.3}Let $(n,\bar Z)$ be a pair 
consisting of an integer $n=0,1,2,\dots,\infty$ and a
spherical operator $\bar Z=(Z_1,\dots,Z_d)$ (which may 
be the nil $d$-tuple when $n\geq 1$).  
Then there is a unique 
nondegenerate representation $\pi$ of $\Cal T_d$ satisfying 
$$
\pi(S_k) = n\cdot S_k\oplus Z_k,\qquad k=1,\dots,d.  
$$
Every nondegenerate representation of $\Cal T_d$ 
on a separable Hilbert space arises in this
way, and if $(n^\prime,\bar Z^\prime)$ is another such pair 
giving rise to a representation $\pi^\prime$, then $\pi^\prime$ is 
unitarily equivalent to $\pi$ iff $n^\prime=n$ and 
$\bar Z^\prime$ is unitarily equivalent to $\bar Z$.  
\endproclaim

\remark{Remarks}
Of course, if $\bar Z$ is the nil $d$-tuple then its corresponding
summand in the definition of $\pi$ is absent.  Let 
$\Cal S\subseteq\Cal B(H)$ be a set of operators acting on a 
Hilbert space $H$.  A subspace $K\subseteq H$ is said to be 
{\it co-invariant} under $\Cal S$ is $\Cal S^*K\subseteq K$.  
$K$ is co-invariant iff its orthogonal 
complement is invariant, $\Cal SK^\perp\subseteq K^\perp$.  
A co-invariant subspace
$K$ is called {\it full}  if $H$ is spanned by 
$\{T\xi\: \xi\in K\}$ where $T$ ranges over the \cstar\ generated
by $\Cal S$.  The following
are equivalent for any co-invariant subace $K$:
\roster
\item"{(8.4.1)}" $K$ is full.  
\item"{(8.4.2)}" $H$ is the smallest reducing subspace for $\Cal S$ 
which contains $K$.  
\item"{(8.4.3)}" For every operator $T$ in the commutant of $\Cal S\cup\Cal S^*$
we have 
$$
TK=\{0\} \implies T=0.
$$  
\endroster
Let $\Cal A$ be the algebra generated by $\Cal S$ and the identity.  
We will often have a situation in which the \cstar\ generated
by $\Cal A$ is spanned by the set of products $\Cal A\Cal A^*$, and
in that case the following criterion can be added to the preceding list
\roster
\item"{(8.4.4)}" H is the smallest invariant subspace for $\Cal S$ which
contains $K$.  
\endroster
Indeed, since $C^*(\Cal A)$ is spanned by $\Cal A\Cal A^*$ 
we have 
$$
\overline{\text{span}}\,C^*(\Cal A)K=
\overline{\text{span}}\,\Cal A\Cal A^*K = \overline{\text{span}}\,\Cal AK,
$$
hence (8.4.1) and (8.4.4) are equivalent.  

Since the $d$-shift is a $d$-contraction, any $d$-tuple 
$(T_1,\dots,T_d)$ of the form 
$$
T_k = n\cdot S_k\oplus Z_k
$$
described in Proposition 8.3 is a $d$-contraction.  
If $K$ is any co-invariant subspace for $\{T_1,\dots,T_d\}$ then 
the $d$-tuple $(T_1^\prime,\dots,T_d^\prime)$
obtained by compressing to $K$
$$
T_k^\prime = P_KT_k\restriction_K
$$
is also a $d$-contraction.  Indeed, 
for each $k=1,\dots,d$ we have 
$$
T^\prime_kT^{\prime *}_k = P_KT_kP_KT_k^*\restriction_K
\leq P_KT_kT_k^*\restriction_K, 
$$
and therefore $\sum_kT_k^\prime T_k^{\prime *}\leq \bold 1$.  
The following implies that 
$d$-tuples obtained from this construction are 
the most general $d$-contractions.  
\endremark

\proclaim{Theorem 8.5}Let $d=1,2,\dots$, let 
$\bar T=(T_1,\dots,T_d)$ be a $d$-contraction
acting on a separable Hilbert space and    
let $\bar S=(S_1,\dots,S_d)$ be the $d$-shift.  Then there is a 
triple $(n,\bar Z,K)$ consisting of an integer $n=0,1,2,\dots,\infty$,  
a spherical operator $\bar Z$, and a full co-invariant subspace
$K$ for the operator
$$
n\cdot \bar S\oplus \bar Z
$$
such that $\bar T$ is unitarily equivalent to the compression 
of $n\cdot \bar S\oplus \bar Z$ to $K$.  

Let $T^\prime=(T_1^\prime,\dots,T_d^\prime)$ be another $d$-contraction
associated with another such triple $(n^\prime,\bar Z^\prime,K^\prime)$.
If $\bar T$ and $\bar T^\prime$ are unitarily equivalent then 
$n^\prime=n$, and there are unitary operators $V\in\Cal B(n\cdot H^2)$ 
and $W:H_{\bar Z}\to H_{\bar Z^\prime}$ such that for $k=1,\dots,d$ 
we have 
$$
VS_k=S_kV, \qquad WZ_k=Z_k^\prime W, 
$$
and which relate $K$ to $K^\prime$ by way of $(V\oplus W)K=K^\prime$.  

Finally, the integer $n$ is the rank of the defect operator 
$$
\bold 1-T_1T_1^*-\dots-T_dT_d^*,
$$ 
and $\bar Z$ is the nil spherical operator iff $\bar T$ is a 
null $d$-contraction.  
\endproclaim

\remark{Remark 8.6}
Notice that the situation of (8.4.4) prevails in this case, and 
we may conclude that for the triple $(n,\bar Z,K)$ associated with 
$\bar T$ by Theorem 8.5, the Hilbert space $\tilde H$ 
on which $n\cdot \bar S\oplus \bar Z$ acts 
is generated as follows
$$
\tilde H = \overline{\text{span}}
\{f(n\cdot S_1\oplus Z_1,\dots,n\cdot S_d\oplus Z_d)\xi\:\xi\in K,
f\in \Cal P\}
$$
$\Cal P$ denoting the set of all polynomials in $d$ complex variables.   
\endremark

Befire giving the proof of Theorem 8.5 we want to emphasize
the following general observation which asserts that, under
certain conditions, a unitary operator which intertwines 
two representations of a subalgebra $\Cal A$ of a 
\cstar\ $\Cal B$ can be extended to a unitary operator 
which intertwines $*$-representations of $\Cal B$.  

We recall a general theorem of Stinespring, which asserts 
that every completely positive map 
$$
\phi: B\to\Cal B(H)
$$
defined on a unital \cstar\ $B$ can be represented in 
the form $\phi(x)=V^*\pi(x)V$, where $\pi$ is a 
representation of $B$ on another Hilbert space 
$H_\pi$, and $V\in\Cal B(H,H_\pi)$.  The pair 
$(V,\pi)$ is called {\it minimal} if 
$$
H_\pi = \overline{\text{span}}[\pi(x)\xi: x\in B,\xi\in H].  
$$
One can always arrange that $(V,\pi)$ is minimal by 
cutting down to a suitable subrepresentation of $\pi$.  

\proclaim{Lemma 8.6}
Let $B$ be a \cstar\ and let $A$ be a  (perhaps
non self-adjoint) subalgebra of $B$ such that 
$$
B = \overline{\text{span}}^{\|\cdot\|}A A^*.  \tag{8.7}
$$
For $k=1,2$ let $\phi_k:B\to \Cal B(H_k)$ be 
$A$-morphisms, and let $U: H_1\to H_2$ be a unitary 
operator such that 
$$
U\phi_1(a) = \phi_2(a)U, \qquad a\in A.  
$$

Let $(V_k,\pi_k)$ be a minimal Stinespring pair for 
$\phi_k$, $\phi_k(x)=V_k^*\pi_k(x)V_k$, $x\in B$.  Then 
there is a unique unitary operator $W: H_{\pi_1}\to H_{\pi_2}$ 
such that 
\roster
\item"{(i)}"
$W\pi_1(x)=\pi_2(x)W, \qquad x\in B$ and 
\item"{(ii)}"
$WV_1 = V_2U$.  
\endroster
\endproclaim

\demo{proof}
Since both $\phi_1$ and $\phi_2$ are $\Cal A$-morphisms, 
the hypothesis on $U$ implies that $U\phi_1(ab^*)=\phi_2(ab^*)U$
for all $a,b\in\Cal A$.  Hence (8.7) implies that 
$U\phi_1(x)=\phi_2(x)U$ for every $x\in B$.  The 
rest now follows from standard uniqueness assertions about
minimal completely positive dilations of completely positive 
maps of \cstar s \cite{3}\qed
\enddemo

\remark{Remark}
There are many examples of subalgebras $A$ of \cstar s $B$ 
that satisfy (8.7) besides the algebra $\Cal A$ of polynomials
in the Toeplitz algebra $\Cal T_d$.  
Indeed, if $\Cal A$ is any algebra of 
operators on a Hilbert space which satisfies 
$$
\Cal A^*\Cal A\subseteq \Cal A + \Cal A^*
$$
then the linear span of $\Cal A\Cal A^*$ is closed under 
multiplication, and hence the norm-closed linear span of 
$\Cal A\Cal A^*$ is a \cstar.  Such examples arise in the 
theory of $E_0$-semigroups \cite{6}, and in the Cuntz 
\cstar s $\Cal O_n$, $n=2,\dots,\infty$.  
\endremark

\demo{proof of Theorem 8.5}
Suppose that the operators $T_k$ act on a Hilbert space $H$.  
Let $\Cal A$ be the algebra of all polynomials in the 
$d$-shift $\bar S=(S_1,\dots,S_d)$.  
By Theorem 6.3 there is an $\Cal A$-morphism 
$$
\phi: \Cal T_d\to \Cal B(H)
$$
such that $\phi(S_k) = T_k$ for $k=1,\dots,d$.  Let 
$$
\phi(X) = V^*\pi(X)V, \qquad X\in \Cal T_d
$$
be a minimal Stinespring representation of $\phi$.  We have 
$$
V^*V=V^*\pi(\bold 1)V = \phi(\bold 1) = \bold 1,
$$ 
hence $V$ is an isometry.  

We claim that $VH$ is co-invariant under $\pi(\Cal A)$, 
$$
\pi(\Cal A)^*VH \subseteq VH.  \tag{8.8}
$$
Indeed, if $A\in\Cal A$ and $P$ denotes the projection 
$P=VV^*$ then for every $X\in \Cal T_d$ we have 
$$
P\pi(A)P\pi(X)V = V\phi(A)\phi(X)=V\phi(AX) 
= P\pi(AX)V=P\pi(A)\pi(X)V
$$
and hence the operator $P\pi(A)P-P\pi(A)$ vanishes on 
$$
\overline{\text{span}}[\pi(X)\xi: X\in \Cal T_d,\xi\in H]=H_{\pi}.  
$$
Thus $\pi(A)^*P = P\pi(A)^*P$ and (8.8) follows.  

Because of minimality of $(V,\pi)$ it follows that 
the subspace $K=VH\subseteq H_{\pi}$ is a {\it full}
co-invariant subspace for the operator algebra $\pi(\Cal A)$.  

Proposition 8.3 shows that if we replace $\pi$ with 
a unitarily equivalent representation and adjust $V$ 
accordingly then we may assume 
that there is an integer 
$n=0,1,2,\dots,\infty$ and a (perhaps nil) spherical operator 
$\bar Z=(Z_1,\dots,Z_d)$ such that $H_\pi=n\cdot H^2\oplus H_{\bar Z}$ 
and 
$$
\pi(S_k) = n\cdot S_k\oplus Z_k, \qquad k=1,\dots,d.  
$$
That proves the first paragraph of Theorem 8.5.  

The second paragraph follows after a straightforward application 
of Lemma 8.6, once one notes that if we are given two triples
$(n,\bar Z,K)$ and $(n^\prime,\bar Z^\prime,K^\prime)$ and we 
define representations $\pi$ and $\pi^\prime$ of $\Cal T_d$ by
$$
\align
\pi(S_k) &= n\cdot S_k\oplus Z_k=\sigma_1(S_k)\oplus\sigma_2(S_k), \\
\pi^\prime(S_k) &= n^\prime\cdot S_k\oplus Z_k^\prime =
\sigma_1^\prime(S_k)\oplus\sigma_2^\prime(S_k), 
\endalign
$$
then $\sigma_1$ is disjoint from $\sigma_2$, 
$\sigma_1^\prime$ is disjoint from $\sigma_2^\prime$, while $\sigma_k$ is 
quasi-equivalent to $\sigma_k^\prime$.  Thus, any unitary operator 
$W$ which intertwines the representations $\pi$ and $\pi^\prime$ must 
decompose into a direct sum $W=W_1\oplus W_2$ where $W_1$ intertwines 
$\sigma_1$ and $\sigma_1^\prime$ and $W_2$ intertwines $\sigma_2$ 
and $\sigma_2^\prime$.  

To prove the third paragraph, choose an integer 
$n=0,1,2,\dots,\infty$, let $\bar Z=(Z_1,\dots,Z_d)$ be 
a spherical operator whose component
operators act on a Hilbert space 
$L$, and let $K\subseteq n\cdot H^2\oplus L$ 
be a full co-invariant subspace for the operator 
$$
n\cdot \bar S\oplus \bar Z, 
$$
where $\bar S=(S_1,\dots,S_d)$ is the $d$-shift.  Define 
$\bar T=(T_1,\dots,T_d)$ by 
$$
T_j=P_K(n\cdot S_j\oplus Z_j)\restriction_K,
$$
$j=1,\dots,d$.  We have to identify the multiplicity $n$ and the 
existence of the spherical summand $\bar Z$ in terms of $\bar T$.  

Let $P_K\in\Cal B(n\cdot H^2\oplus L)$ denote the projection 
on $K$.  Since $K$ is co-invariant under 
$n\cdot \bar S\oplus \bar Z$ we 
have 
$$
P_K( n\cdot S_j\oplus Z_j)=P_K( n\cdot S_j\oplus Z_j)P_K=
T_jP_K
$$
for every $j=1,\dots,d$, and hence 
$$
T_jT_j^* = P_K(n\cdot S_jS_j^*\oplus Z_jZ_j^*)\restriction_K. \tag{8.9}
$$
By the remarks following Definition 2.10 we may sum on $j$
to obtain 
$$
\sum_{j=1}^d T_jT_j^* = 
P_K(n\cdot (\bold 1-E_0)\oplus \bold 1_L)\restriction_K=
\bold 1_K-P_K(n\cdot E_0\oplus 0)\restriction_K \tag{8.10}
$$
where $E_0\in\Cal B(H^2)$ denotes the one-dimensional projection
onto the constants.  

From (8.10) we find that the defect operator $D$ has the form 
$$
D=\bold 1_K-T_1T_1^*-\dots-T_dT_d^* = 
P_K(n\cdot E_0\oplus 0)\restriction_K.  \tag{8.11}
$$
Now for any positive operator $B$ we have $B\xi=0$ iff 
$\<B\xi,\xi\> = 0$.  Thus the relation (8.11) between 
the positive operators $D$ and  $n\cdot E_0\oplus 0$ implies 
that their kernels are related by 
$$
\{\xi\in K: D\xi=0\} = \{\xi\in K: (n\cdot E_0\oplus 0)\xi=0\},   
$$
and hence 
$$
\text{rank} D = \dim ((n\cdot E_0\oplus 0 )K).  
$$  
The dimension of the space $N=(n\cdot E_0\oplus 0 )K$ is easily 
seen to be $n$.  Indeed, notice that if $A\in\Cal B(H^2)$ 
is a polynomial in the operators $S_1,\dots,S_d$ then we have 
$E_0A = E_0AE_0 = \<A1,1\>E_0$, and hence $E_0A$ is a scalar 
multiple of $E_0$.  Similarly, if $B\in \Cal B(n\cdot H^2\oplus L)$
is a polynomial in the operators 
$n\cdot S_1\oplus Z_1,\dots, n\cdot S_d\oplus Z_d$ then 
$(n\cdot E_0\oplus 0)B$ is a scalar multiple of $(n\cdot E_0\oplus 0)$, 
and hence for all such $B$ we have 
$$
(n\cdot E_0\oplus 0)BK \subseteq N.  
$$
Because $K$ is a full co-invariant subspace,  (8.4.4)
implies that $n\cdot H^2\oplus L$ is spanned by vectors of the 
form $B\xi$, with $B$ as above and $\xi\in K$.  It follows that 
$$
(n\cdot E_0\oplus 0 )(n\cdot H^2\oplus L)\subseteq N,
$$
and therefore $N$ is the range of the $n$-dimensional projection 
$n\cdot E_0\oplus 0$.  Hence $\dim N=n$.  

Finally, we consider case in which $\bar T$ is a null $d$-contraction. 
Let $Q$ and $P$ be the completely positive maps on $\Cal B(H^2)$ and
$\Cal B(K)$ given respectively by 
$$
\align
P(A)&=S_1AS_1^*+\dots+ S_dAS_d^*, \qquad  A\in \Cal B(H^2),\\
Q(B)&=T_1BT_1^*+\dots+T_dBT_d^*, \qquad B\in \Cal B(K).
\endalign  
$$ 
Formula (8.9) implies that
$Q(\bold 1_K) = 
P_K (n\cdot P(\bold 1_{H^2})\oplus \bold 1_L)\restriction_K$.  
Similarly, using co-invariance of $K$ repeatedly 
as in (8.9) we have 
$$
T_{j_1}\dots T_{j_r}T_{j_r}^*\dots T_{j_1}^* = 
P_K(n\cdot( S_{j_1}\dots S_{j_r}S_{j_r}^*\cdots S_{j_1}^*)\oplus 
Z_{j_1}\cdots Z_{j_r}Z_{j_r}^*\cdots Z_{j_1}^*)\restriction_K
$$
for every $j_1,\dots,j_r\in\{1,\dots,d\}$.  After summing 
on $j_1,\dots,j_r$ we obtain 
$$
Q^r(\bold 1_K) = P_K(n\cdot P^r(\bold 1_{H^2})\oplus 
\bold 1_L)\restriction_K,\qquad r=1,2,\dots.  
$$
Since $P^r(\bold 1_{H^2})\downarrow 0$ as $r\to\infty$, we have 
$$
\lim_{r\to\infty}Q^r(\bold 1_K) = P_K(0\oplus \bold 1_L)\restriction_K.  
$$
We conclude that $\bar T$ is a null $d$-tuple iff $0\oplus L\perp K$, 
that is, $K\subseteq n\cdot H^2\oplus \{0\}$.  
Noting that $n\cdot H^2\oplus \{0\}$ is a reducing subspace 
for the operator $n\cdot \bar S\oplus\bar Z$ we see from 
(8.4.2) that 
$$
n\cdot H^2\oplus L \subseteq n\cdot H^2\oplus \{0\},   
$$
and therefore $L=\{0\}$.  But a 
spherical $d$-tuple cannot be the zero $d$-tuple except when it 
is the nil $d$-tuple, and thus we have proved that $\bar T$ is
a null $d$-contraction iff $\bar Z$ is nil\qed
\enddemo

The two extreme cases of Theorem 8.5 in which $n=0$ and $n=1$ 
are noteworthy.  From the case $n=0$ we deduce the following 
result of Athavale \cite{11}, which was established by  
entirely different methods.  

\proclaim{Corollary 1}
Let $T_1,\dots,T_d$ be a set of commuting operators on a Hilbert 
space $H$ such that $T_1^*T_1+\dots+T_d^*T_d=\bold 1$.  Then
$(T_1,\dots,T_d)$ is a subnormal $d$-tuple.  
\endproclaim
\demo{proof}
Let $A_k=T_k^*$.  $(A_1,\dots,A_d)$ is a $d$-contraction for which 
$$
n=\text{rank}(\bold 1-A_1A_1^*-\dots-A_dA_d^*)=0.
$$  
Theorem 8.5 implies that there is a spherical operator 
$\bar Z=(Z_1,\dots,Z_d)$ acting on a Hilbert space 
$\tilde H\supseteq H$ such that $Z^*_kH\subset H$ and 
$A_k$ is the compression of $Z_k$ to $H$, $k=1,\dots,d$.  
Hence $T_k=A_k^*=Z_k^*\restriction_H$ for every $k$, so 
that $(Z_1^*,\dots,Z_d^*)$ is a normal $d$-tuple which 
extends $(T_1,\dots,T_d)$ to a larger Hilbert space\qed
\enddemo

From the case $n=1$ we have the following description  
of all $d$-contractions that can be 
obtained by compressing the $d$-shift
to a co-invariant subspace.  

\proclaim{Corollary 2}
Every nonzero co-invariant subspace $K\subseteq H^2$ for the 
$d$-shift $\bar S=(S_1,\dots,S_d)$ is full, and the compression
of $\bar S$ to $K$
$$
T_k = P_KS_k\restriction_K, \qquad k=1,\dots,d
$$
defines a null $d$-contraction $\bar T=(T_1,\dots,T_d)$ for 
which 
$$
\text{\rm{rank}}(\bold 1-T_1T_1^*-\dots-T_dT_d^*)=1.\tag{8.12}
$$  
If $K^\prime$ is another co-invariant subspace for $\bar S$ 
which gives rise to $\bar T^\prime$, then $\bar T$ and 
$\bar T^\prime$ are unitarily equivalent if and only if 
$K=K^\prime$.  

Every null $d$-contraction $(T_1,\dots,T_d)$ satisfying (8.12) 
is unitarily equivalent to one obtained by compressing 
$(S_1,\dots,S_d)$ to a co-invariant subspace of $H^2$.  
\endproclaim

\demo{proof}
Let $\{0\}\neq K\subseteq H^2$ be a 
a co-invariant subspace for the set 
of operators $\{S_1,\dots,S_d\}$.  
Since $\Cal T_d$ is an irreducible \cstar\ it follows 
that $K$ satisfies 
condition (8.4.2), hence it is full.  Let $T_j$ be the compression 
of $S_j$ to $K$, $j=1,\dots,d$.  
The canonical triple associated with $\bar T=(T_1,\dots,T_d)$ is therefore
$(1,\text{nil},K)$, and the third paragraph of Theorem 8.5 
implies that $\bar T$ is a null $d$-contraction satisfying (8.12).

If $K^\prime$ is another co-invariant subspace of $H^2$ 
giving rise to a $d$-contraction $\bar T^\prime$ which is 
unitarily equivalent to $\bar T$ then Theorem 8.5 implies 
that there is a unitary operator $V$ which commutes with 
$\Cal S=\{S_1,\dots,S_d\}$ such that $VK=K^\prime$.  
Because $V$ is unitary it must commute with $\Cal S^*$ as 
well, and hence with the Toeplitz algebra $\Cal T_d$.  The 
latter is irreducible, hence $V$ must be a scalar multiple 
of the identity operator, hence $K^\prime=K$.  

Finally, if $\bar T=(T_1,\dots,T_d)$ is any null $d$-contraction
then Theorem 8.5 implies that 
the spherical summand $\bar Z$ of its dilation must 
be the nil $d$-tuple, and if in addition 
$$
\text{rank}(\bold 1-T_1T_1^*-\dots-T_dT_d^*)=1,
$$
then the canonical triple associated with $\bar T$ is 
$(1,\text{nil},K)$ for some subspace $K$ 
of $H^2$ which is co-invariant under the $d$-shift \qed
\enddemo

Lemma 7.13 asserts that the identity representation of the Toeplitz
\cstar\ is a boundary representation for the unital operator 
space generated by the $d$-shift.  This fact has a number 
of interesting consequences, and we conclude with 
a brief discussion of one of them.  
Rudin posed the following function-theoretic problem in the sixties: 
Do there exist nonconstant inner functions in 
$H^\infty(B_d)$ \cite{37}?  This problem was 
finally solved (affirmatively) in 1982 by B. A. Aleksandrov 
\cite{38}.  
The following Proposition implies 
that the answer to the analogue of Rudin's 
question for the multiplier algebra $\Cal M$ is the opposite: 
there are no nontrivial isometries in $\Cal B(H^2)$
which commute with $\{S_1,\dots,S_d\}$ when $d\geq 2$.  
Indeed, we have the following more general assertion.  

\proclaim{Proposition 8.13}
Let $T_1,T_2,\dots$ be a finite or infinite 
sequence of operators on $H_d^2$,
$d\geq 2$, which 
commute with the $d$-shift and which satisfy 
$$
T_1^*T_1+T_2^*T_2+\dots = \bold 1.  \tag{8.6}
$$
Then each $T_j$ is a scalar multiple of the identity operator.  
\endproclaim  
\demo{proof}
Consider the completely positive linear map $\phi$ defined 
on $\Cal B(H^2)$ by 
$$
\phi(A)=T_1^*AT_1+T_2^*AT_2+\dots.  
$$
The sum converges strongly for every operator $A$ because by 
(8.6) we have 
$$
\|T_1\xi\|^2+\|T_2\xi\|^2+\dots = \|\xi\|^2<\infty,\qquad \xi\in H^2.  
$$
Moreover, since each $T_k$ commutes with each $S_j$ we 
have $T_k^*S_jT_k = T_k^*T_kS_j$, and thus from (8.6) we 
conclude that $\phi(A)=A$ for every $A$ in 
$\Cal S=\text{span}\{\bold 1,S_1,\dots,S_d\}$.  Since the 
identity representation of $\Cal T_d$ is a boundary representation
for $\Cal S$ it follows that $\phi(A)=A$ for every $A$ in 
the Toeplitz \cstar\ $\Cal T_d$.  

Let $n$ be the number of operators in the 
sequence $T_1,T_2,\dots$ and  let $V$ be the 
linear map of $H^2$ to $n\cdot H^2$ defined by 
$$
V\xi = (T_1\xi,T_2\xi,\dots).  
$$
Because of (8.6) $V$ is an isometry.  Letting $\pi$ be the 
representation of $\Cal B(H^2)$ on $n\cdot H^2$ defined by 
$$
\pi(A)=A\oplus A\oplus \dots, 
$$
we find that $(V,\pi)$ is a Stinespring pair for $\phi$, 
$$
\phi(A) = V^*\pi(A)V, \qquad A\in \Cal B(H^2).  
$$
Since 
$$
(VA-\pi(A)V)^*(VA-\pi(A)V) = 
A^*\phi(\bold 1)A -\phi(A)^*A-A^*\phi(A)+\phi(A^*A)
= 0, 
$$
we conclude that $VA-\pi(A)V=0$.  
By examining the components of this 
operator equation one sees that $T_kA= AT_k$ for every $k$ and 
every $A\in\Cal T_d$.  Since $\Cal T_d$ is an irreducible 
\cstar\ it follows that each $T_k$ must be a scalar 
multiple of the identity operator\qed
\enddemo


\subheading{Appendix. Trace estimates}

Fix $d=1,2,\dots$,  let $E_d$ be a $d$-dimensional Hilbert 
space, and let 
$$
\Cal F_+(E_d)=\Bbb C\oplus E_d\oplus E_d^2\oplus\dots
$$
be the symmetric Fock space over $E_d$.  The number operator is 
the unbounded self-adjoint diagonal operator $N$ satisfying 
$N\xi = n\xi$, $\xi\in E_d^n$, $n=0,1,\dots$.  Let $P_n$ be the 
projection on $E_d^n$.  Then for every $p>0$,  $(\bold 1+N)^{-p}$ 
is a positive compact operator
$$
(\bold 1+N)^{-p} = \sum_{n=0}^\infty (n+1)^{-p}P_n 
$$
whose trace is given by
$$
\text{trace}(\bold 1+N)^{-p} = 
\sum_{n=0}^\infty \frac{\dim E_d^n}{(n+1)^p}.  \tag{A.1}
$$
Thus $(\bold 1+N)^{-1}$ belongs to the Schatten class 
$\Cal L^p(\Cal F_+(E_d))$ iff the infinite series (A.1) converges.  
In this appendix we show that is the case iff $p>d$.  
Notice that the function of a complex variable defined for $\Re z>d$ by
$$
\zeta_d(z) = \text{trace}(\bold 1+N)^{-z}
$$
is a $d$-dimensional variant of the Riemann zeta function, since
for $d=1$ we have $\dim E_d^n=1$ for all $n$ and hence 
$$
\zeta_1(z) = \sum_{n=1}^\infty \frac{1}{n^z}.  
$$

We calculate the generating function for the 
coefficients $\dim E_d^n $.  

\proclaim{Lemma A.2}
The numbers $a_{n,d}= \dim E_d^n$ are the coefficients of the 
series expansion 
$$
(1-z)^{-d} = \sum_{n=0}^\infty a_{n,d}z^n, \qquad |z|<1.  
$$
\endproclaim
\demo{proof}
Note that the numbers $a_{n,d}$ satisfy the 
recurrence relation
$$
a_{n,d+1} = a_{0,d}+a_{1,d}+\dots+a_{n,d},
\qquad n=0,1,\dots, d=1,2,\dots.  \tag{A.3}  
$$
Indeed, if we choose a basis $e_1,\dots,e_d$ for $E_d$ then 
the set of symmetric products
$$
\{e_{i_1}e_{i_2}\dots e_{i_n}: 1\leq i_1\leq \dots \leq i_n\leq d\}
$$
forms a basis for the vector space $E_d^n$ and hence 
$a_{n,d}$ is the cardinality of the set 
$$
S_{n,d}=\{(i_1,\dots,i_n)\in\{1,\dots,d\}^n : 
1\leq i_1\leq\dots\leq i_n\leq d\}.  
$$
Since $S_{n,d+1}$ decomposes into a disjoint union 
$$
S_{n,d+1} = \bigsqcup_{k=0}^n\{(i_1,\dots,i_n)\in S_{n,d+1}: 
i_k\leq d,i_{k+1}=\dots=i_n=d+1 \}, 
$$
and since the $k$th set on the right has the same cardinality $a_{k,d}$
as $S_{k,d}$, (A.3) follows.  

From (A.3) we find that $a_{n,d+1}-a_{n-1,d+1}=a_{n,d}$.  
Thus if we let $f_d$ be the formal power series
$$
f_d(z) = \sum_{n=0}^\infty a_{n,d}z^n  \tag{A.4}
$$
then $f_{d+1}(z) - zf_{d+1}(z)=f_d(z)$, hence
$$
f_{d+1}(z) = \frac{f_d(z)}{1-z}.  
$$
Lemma A.1 follows after noting that 
$f_1(z)=1+z+z^2+\dots=(1-z)^{-1}$\qed
\enddemo

\remark{Remark}
Notice that the power series of Lemma A.2
converges absolutely 
to the generating function $(1-z)^{-d}$ throughout
the open unit disk $|z|<1$.
\endremark

By evaluating successive 
derivatives of the generating function at the origin, 
we find that 
$$
\dim E_d^n = \frac{(n+d-1)(n+d-2)\dots d}{n!}=
\frac{(n+d-1)!}{n!(d-1)!}. \tag{A.5}
$$
A straightforward application of Stirling's 
formula \cite{36, p. 194} 
$$
N!\sim \sqrt{2\pi}\,N^{N+1/2}e^{-N}
$$
leads to 
$$
\lim_{n\to\infty}(n+1)^{-d+1}\frac{(n+d-1)!}{n!} = 
\frac{1}{(d-1)!}
$$
and hence
$$
\dim E_d^n \sim \frac{(n+1)^{d-1}}{(d-1)!}.  \tag{A.6}
$$

We now prove the assertion of (5.2).  
\proclaim{Theorem}
For $p>0$ we have $\text{\rm{trace}}(\bold 1+N)^{-p}<\infty$ if and 
only if $p>d$.  
\endproclaim
\demo{proof}
By (A.6), the infinite series 
$$
\text{trace}(\bold 1+N)^{-p} = \sum_{n=0}^\infty\frac{\dim E_d^n}{(n+1)^p}
$$
converges if and only if the series 
$$
\sum_{n=0}^\infty \frac{1}{(n+1)^{p-d+1}}
$$
converges; i.e., iff $p>d$\qed
\enddemo

\vfill
\pagebreak

\end
\Refs
%


\ref\no 1\by Agler, Jim \paper The Arveson extension theorem 
and coanalytic models\jour  Int. Eq. and Op. Th.
\vol 5 \yr 1982\pages 608--631
\endref

\ref\no 2\bysame \paper Hypercontractions and
subnormality\jour  Jour. Op. Th.\vol 13 \yr 1985
\pages 203--217
\endref



\ref\no 3\by Arveson, W.\paper Subalgebras of $C^*$-algebras
\jour Acta Math\vol 123 \yr 1969\pages 141--224
\endref

\ref\no 4\bysame, W.\paper Subalgebras of $C^*$-algebras II
\jour Acta Math\vol 128\yr 1972\pages 271--308
\endref

\ref\no 5\bysame\book An Invitation to $C^*$-algebras
\bookinfo Gratuate Texts in Mathematics \publ
Springer-Verlag \vol 39\year 1976
\endref


\ref\no 6\bysame\paper \cstar s associated with sets
of semigroups of operators
\yr 1991\jour Int. J. Math. \vol 2, no. 3\pages 235--255
\endref



\ref\no 7\bysame\paper Minimal $E_0$-semigroups\inbook 
Operator Algebras and their Applications\publ AMS
\bookinfo Fields Institute Communications \ed Fillmore, 
P. and Mingo, J.  \yr 1997\pages 1--12
\endref

\ref\no 8\bysame\paper The index of a quantum dynamical
semigroup\jour J. Funct. Anal. \paperinfo to appear
\endref

\ref\no 9\bysame\paper On the index and dilations 
of completely positive semigroups\jour Int. J. Math.
\paperinfo to appear
\endref

\ref\no 10\bysame\paper Pure $E_0$-semigroups and
absorbing states\jour Comm. Math. Phys.  
\paperinfo to appear
\endref


\ref\no 11\by Athavale, A.\paper On the intertwining
of joint isometries
\jour Jour. Op. Th. \vol 23\yr 1990\pages 339--350
\endref

\ref\no 12\bysame\paper Model theory on 
the unit ball of $\Bbb C^n$
\jour Jour. Op. Th. \vol 27\yr 1992
\pages 347--358
\endref

\ref\no 13\by Attele, K. R. M. and Lubin, A. R. 
\paper Dilations and commutant lifting for jointly
isometric operators--a geometric approach
\jour J. Funct. Anal.\yr 1996 \vol 140 
\pages 300--311
\endref



\ref\no 14\by Bhat, B. V. R. \paper An index theory for quantum 
dynamical semigroups\jour Trans. AMS \vol 348 \yr 1996 
\pages 561-583
\endref


\ref\no 15\by Bunce, John W.\paper Models for $n$-tuples 
of noncommuting operators
\jour J. Funct. Anal. \yr 1984\vol 57\pages 21--30
\endref


\ref\no 16\by Coburn, L. A. \paper Singular integral operators 
and Toeplitz operators on odd spheres
\jour Ind. Univ. Math. J.\yr 1973 \vol 23\pages 433--439
\endref

\ref\no 17\by Curto, R. and Vasilescu, F.-H.
\paper Automorphism invariance of the operator-valued
Poisson transform 
\jour Acta Sci. Math. (Szeged)\yr 1993 \vol57\pages 65--78
\endref

\ref\no 18\by Davidson, K. R. and Pitts, D.
\paper Invariant subspaces and hyper-reflexivity
for free semigroup algebras\paperinfo preprint
\endref

\ref\no 19\bysame
\paper The algebraic structure of non-commutative
analytic Toeplitz algebras \paperinfo preprint
\endref


\ref\no 20\by Drury, S.\paper A generalization of 
von Neumann's inequality to the complex ball
\jour Proc. AMS\yr 1978 \vol 68\pages 300--304
\endref


\ref\no 21\by Frazho, Arthur E.\paper Models for 
noncommuting operators
\jour J. Funct. Anal. \yr 1982 \vol 48\pages 1--11
\endref

\ref\no 22\by Halmos, P. R.\book A Hilbert space
problem book
\publ Van Nostrand
\yr 1967 \publaddr New York
\endref

\ref\no 23\by Hamana, M.\paper Injective envelopes
of $C^*$-algebras\jour J. Math. Soc. Japan
\yr 1979 \vol 31\pages 181--197
\endref

\ref\no 24\bysame\paper Injective envelopes 
of operator systems
\jour Publ. R.I.M.S., Kyoto Univ.\yr 1979
\vol 15, no. 3\pages 773--785
\endref


\ref\no 25\by M\"uller, V. and Vasilescu, F.-H.
\paper Standard models for some commuting multioperators
\jour Proc. AMS\yr 1993 \vol 117\pages 979--989
\endref


\ref\no 26\by Paulsen, V.\book Completely bounded maps 
and dilations\publ Wiley\yr 1986\publaddr New York
\endref


\ref\no 27\by Pisier, G.\book Similarity problems
and completely bounded maps
\publ Springer Verlag Lecture Notes in Mathematics\vol 
1618\yr 1995
\endref

\ref\no 28\by Popescu, G.\paper Models for infinite sequences 
of noncommuting operators
\jour Acta Sci. Math. (Szeged)\yr 1989 \vol 53\pages 355--368
\endref

\ref\no 29\bysame\paper Isometric dilations for 
infinite sequences of noncommuting operators
\jour Trans. Amer. Math. Soc.\yr 1989 \vol 316\pages 523--536
\endref

\ref\no 30\bysame\paper von Neumann inequality
for $(\Cal B(H)^n)_1$
\jour Math. Scand.\yr 1991 \vol 68\pages 292--304
\endref

\ref\no 31\bysame\paper On intertwining dilations 
for sequences of noncommuting operators
\jour J. Math. Anal. Appl.\yr 1992 \vol 167\pages 382--402
\endref


\ref\no 32\bysame\paper Functional calculus for 
noncommuting operators
\jour Mich. Math. J.\yr 1995 \vol 42\pages 345--356
\endref

\ref\no 33\bysame\paper Multi-analytic operators on 
Fock space \jour Math. Ann.\yr 1995 \vol 303\pages 31--46
\endref


\ref\no 34\bysame\paper Noncommutative disc algebras 
and their representations
\jour Proc. Amer. Math. Soc.\yr 1996 \vol 124\pages 2137--2148
\endref

\ref\no 35\bysame\paper Poisson transforms on some 
$C^*$-algebras generated by isometries
\yr 1995\paperinfo preprint
\endref


\ref\no 36\by Rudin, W.\book Principles of Mathematical Analysis,
3rd edition\publ McGraw Hill\yr 1976
\endref


\ref\no 37\bysame\book Function theory in the unit ball of 
$\Bbb C^n$\publ Springer Verlag\yr 1980
\endref

\ref\no 38\bysame\book New constructions of functions holomorphic
in the unit ball of $\Bbb C^n$\bookinfo CBMS publication number 63
\publ AMS\yr 1986
\endref


\ref\no 39\by Sarason, D. \paper On spectral sets having
connected complement
\jour Acta Sci. Math.\vol 26 \yr 1965\pages 289--299
\endref

\ref\no 40\by Segal, I. E. \paper Tensor algebras over 
Hilbert spaces\jour Trans. Amer. Math. Soc. 
\yr 1956\pages 106--134
\endref

\ref\no 41\bysame\paper Tensor algebras over Hilbert 
spaces II\jour Annals Math.\vol 63, no. 1 \yr 1956\pages 160--175
\endref

\ref\no 42\by SeLegue, Dylan \paperinfo Berkeley Ph.D. 
dissertation\yr 1997.  
\endref

\ref\no 43\by Sz.-Nagy, B. and Foias, C.\book
Harmonic analysis of operators on Hilbert space 
\publ American Elsevier, New York\yr 1970
\endref

\ref\no 44\by Vasilescu, F.-H.\paper An operator-valued
Poisson kernel 
\jour Jour. Funct. Anal.\vol 110\yr 1992\pages 47--72
\endref


\ref\no 45\bysame\paper Operator-valued 
Poisson kernels and standard models in several variables
\inbook Algebraic methods in operator theory\ed Ra\'ul 
Curto and Palle Jorgensen\yr 1994 \publ Birkh\"auser 
\pages 37--46
\endref

\endRefs
